\documentclass[11pt,reqno]{amsart}

\usepackage[T1]{fontenc}
\usepackage[utf8]{inputenc}
\usepackage{lmodern}
\usepackage{microtype}
\usepackage{amsmath,amssymb,amsfonts,mathtools}
\usepackage{mathrsfs}
\usepackage{enumitem}
\usepackage{xcolor}
\usepackage{booktabs}
\usepackage{tikz-cd}
\usepackage[hidelinks]{hyperref}

\makeatletter
\@namedef{subjclassname@2020}{\textup{2020} Mathematics Subject Classification}
\makeatother

\numberwithin{equation}{section}

\newtheorem{theorem}{Theorem}[section]
\newtheorem{proposition}[theorem]{Proposition}
\newtheorem{lemma}[theorem]{Lemma}
\newtheorem{corollary}[theorem]{Corollary}
\newtheorem{conjecture}[theorem]{Conjecture}

\theoremstyle{definition}
\newtheorem{definition}[theorem]{Definition}
\newtheorem{construction}[theorem]{Construction}

\theoremstyle{remark}
\newtheorem{remark}[theorem]{Remark}

\newcommand{\Q}{\mathbf Q}
\newcommand{\Z}{\mathbf Z}

\newcommand{\C}{\mathbf C}

\newcommand{\Hh}{\mathfrak H}
\newcommand{\G}{\mathrm{GL}_2}
\newcommand{\SL}{\mathrm{SL}_2}
\newcommand{\Mat}{\mathrm M_2}
\newcommand{\Div}{\operatorname{Div}}
\newcommand{\Dist}{\operatorname{Dist}}
\newcommand{\Symb}{\operatorname{Symb}}
\newcommand{\Tay}{\operatorname{Tay}}
\newcommand{\Int}{\operatorname{Int}}

\newcommand{\Eis}{\mathrm{Eis}}

\newcommand{\cF}{\mathcal F}
\newcommand{\cH}{\mathcal H}
\newcommand{\cE}{\mathcal E}

\newcommand{\cS}{\mathcal S}
\newcommand{\cC}{\mathcal C}
\newcommand{\scrE}{\mathscr E}
\newcommand{\scrM}{\mathscr M}
\newcommand{\scrN}{\mathscr N}
\newcommand{\scrS}{\mathscr S}
\newcommand{\scrC}{\mathscr C}

\newcommand{\one}{\mathbf 1}

\newcommand{\prcusp}{\operatorname{pr}_{\mathrm{cusp}}}

\title[An Eisenstein--Kronecker modular symbol]{An Eisenstein--Kronecker Modular Symbol and Beilinson--Kato Distributions}
\author{Vicen\c{t}iu Pa\c{s}ol}
\address{Institute of Mathematics ``Simion Stoilow'' of the Romanian Academy,
21 Calea Grivi\c{t}ei, 010702 Bucharest, Romania}
\email{vicentiu.pasol@imar.ro}
\date{}


\hypersetup{
  pdftitle={An Eisenstein--Kronecker Modular Symbol and Beilinson--Kato Distributions},
  pdfauthor={Vicentiu Pasol},
  pdfsubject={Modular symbols and Eisenstein--Kronecker distributions},
  pdfkeywords={modular symbols, Eisenstein--Kronecker series, Siegel units, Beilinson--Kato elements, distributions}
}

\begin{document}

\begin{abstract}
We construct a $\G(\Q)$-modular symbol $\Phi$ from products of
Eisenstein--Kronecker series by holomorphic projection modulo Eisenstein
series.  Its finite-level specializations recover the Borisov--Gunnells
maps after a finite Fourier transform and polynomial duality.
We reconstruct $\Phi$ by applying local explicit reciprocity to the
ordered Siegel units defining the Beilinson--Kato distribution of
Busuioc, Park, Patashnick, and Stevens.  In weight $k\ge4$, rational
torus covariance and the Manin relations determine the entire symbol
in that weight from partial Fourier data for the coefficients of
$X^{k-2}$ and $Y^{k-2}$.
The cyclic Manin relation also has a geometric consequence.
For every odd prime $p$, its weight-two equations, after elimination of
the weight-two coordinates, cut out the modular curve $X(p)$
scheme-theoretically in weight-one Eisenstein coordinates.  At level
one its coefficients give all linear relations among unordered products
of two Eisenstein series.
\end{abstract}

\maketitle
\enlargethispage{4pt}

\section{Introduction}

We construct a modular symbol $\Phi$ from products of
Eisenstein--Kronecker series and relate it to the Beilinson--Kato
modular symbol of Busuioc, Park, Patashnick, and Stevens (BPPS)
\cite{BPPS}.  Both symbols are equivariant under the full group
$\G(\Q)$.  On the cusps $\mathbf P^1(\Q)$ its matrices act by
fractional linear transformations, with scalar matrices as kernel,
so the induced group is $\operatorname{PGL}_2(\Q)$.  Given distinct
rational cusps $r,s$, choose $g\in\G(\Q)$ carrying $(\infty,0)$ to
$(r,s)$.  For either symbol $\Psi$, equivariance gives
\[
 \Psi([r]-[s])=g\cdot\Psi([\infty]-[0]),
\]
where $g$ acts on the coefficient module of $\Psi$.  Thus one value
determines the entire symbol.  Local explicit reciprocity applied to the
Siegel units defining the BPPS distribution gives part of the finite Fourier data of
$\Phi([\infty]-[0])$.  The rational diagonal torus stabilizes
$(\infty,0)$, and its action, together with the Manin relations,
determines the remaining coefficients in every weight.

The construction of $\Phi$ and its finite-level specializations go back
to the 2006 version of this article
\cite[Theorems~0.1--0.2]{Pasol2006}.  A single distribution thus
encodes the Borisov--Gunnells maps \cite{BG1,BG2} as the weight and
level vary.  In the analytic construction the cyclic Manin relation
arises from a three-term identity for Eisenstein--Kronecker series.
Here we also study the equations furnished by that identity.  Their
coefficients span the space of level-one product relations described
below.  At full prime level their weight-two specialization gives
projective equations for $X(p)$ in the coordinates supplied by
weight-one Eisenstein series, for every odd prime $p$.

Put
\[
 \Delta_0=\Div^0(\mathbf P^1(\Q)),\qquad D_\infty=[\infty]-[0].
\]
For a left $\G(\Q)$-module $W$, a $\G(\Q)$-modular symbol with values
in $W$ is a $\G(\Q)$-equivariant homomorphism $\Delta_0\to W$.
Write $\Symb_{\G(\Q)}(W)$ for the group of such symbols.  This is the
convention of \cite[Introduction]{BPPS}, written with left actions.
Our divisor $D_\infty$ is the negative of that used in
\cite[Theorem~1.1.4]{BPPS}.

Let $\scrC_2$ be the space of formal series whose $X^iY^j$-coefficient
is a modular form of weight $i+j+2$, modulo Eisenstein series;
the coefficient modules and their actions are defined in
Section~\ref{sec:coefficients}.  For a rational lattice coset
$V\subset\Q^2$, let $E_V(X)$ be the generating series of its
Eisenstein--Kronecker coefficients, with the omitted-term regularization
at zero.  Define, with columns as the two factors,
\begin{equation}\label{eq:intro-mu}
 \mu(V_1\times V_2)=\cH\bigl(E_{V_1}(X)E_{V_2}(Y)\bigr),
\end{equation}
where $\cH$ is algebraic holomorphic projection followed by passage modulo
Eisenstein series.  The formula extends by linearity to the test-function
space on $\Mat(\Q)$.

\begin{theorem}\label{thm:intro-main}
There is a unique $\G(\Q)$-modular symbol
\[
 \Phi:\Delta_0\longrightarrow\Dist(\Mat(\Q),\scrC_2)
\]
such that $\Phi(D_\infty)=\mu$.
\end{theorem}

The theta addition formula and the heat equation give a quadratic
identity that survives regularization when a torsion parameter is zero.
After holomorphic projection and passage to the Eisenstein quotient,
its cyclic sum vanishes.  Together with the symmetry exchanging the two
factors, this gives the Manin relations
\[
 (1+S)\mu=0,\qquad (1+R+R^2)\mu=0,
 \qquad
 S=\begin{pmatrix}0&-1\\1&0\end{pmatrix},\quad
 R=\begin{pmatrix}0&-1\\1&-1\end{pmatrix}.
\]
For $N>1$, put
\[
 U_1(N)=\begin{pmatrix}0&0\\0&1/N\end{pmatrix}+\Mat(\Z).
\]

\begin{theorem}\label{thm:intro-BG}
For $N>1$ and $k\ge2$, the map
\[
 D\longmapsto\Phi(D)(U_1(N))^{[k-2]}
\]
is a $\Gamma_1(N)$-modular symbol.  Under the cuspidal splitting,
its values lie in
\[
 S_k(\Gamma_1(N))\otimes\operatorname{Sym}^{k-2}(\C^2).
\]
The weight-$k$ Borisov--Gunnells map is recovered from $\mu$ by finite
Fourier transform on all pairs modulo $N$, followed by polynomial duality
and restriction to primitive pairs, as in Section~\ref{sec:BG}.
\end{theorem}

The reconstruction starts from the ordered product of Siegel
distributions at $[0]-[\infty]$, evaluated on tests avoiding zero in
each row.  We use the explicit reciprocity law for Kato's Euler
system \cite{Kato} in the local form proved by Wang
\cite[Theorem~1.4]{Wang}.  After smoothing and taking moments, the
local exponential of the Kummer cup product gives Colmez's Eisenstein
product distribution $z_{\Eis}(k,j)$ \cite[\S1.3.5]{Colmez}, with the
two factors smoothed and the sign convention of \eqref{eq:Wang-reciprocity}.
One factor carries an additive character where the definition of $\mu$
uses a translated lattice.  Finite symplectic Fourier transform in that
factor relates the two products.  To specify the part of $\mu$ that this
calculation determines directly, introduce the following tests.

For $r\in\Q_{>0}$, $N\ge1$, and $a,b,x\in(\Z/N\Z)^2$, put
\[
 \begin{gathered}
 e^{r,N}_{a,b}=[r(a,b)/N+r\Mat(\Z)],\qquad
 \langle x,a\rangle=x_1a_2-x_2a_1,\\
 f^{r,N}_{x,b}=\sum_a e^{2\pi i\langle x,a\rangle/N}e^{r,N}_{a,b}.
 \end{gathered}
\]
The parentheses $(a,b)$ denote a matrix with columns $a,b$.
Let $\mu_k^{(j)}$ be the coefficient of $X^{k-j-1}Y^{j-1}$ in $\mu$.

\begin{samepage}
\begin{theorem}\label{thm:intro-comparison}
On row tests avoiding zero in each factor, take the BPPS distribution
in its defining presentation by ordered pairs of algebraic Siegel units.
Applying Wang's construction to their Kummer classes as in
Section~\ref{sec:BK} determines all values
\[
 \mu_k^{(j)}(f^{r,N}_{x,b}),\qquad
 k\ge2,\quad 1\le j\le k-1,\quad j\ne2,\quad
 \langle x,b\rangle\not\equiv0\pmod N.
\]
Among distributions with the action \eqref{eq:old-action} satisfying the
torus and Manin relations, these values determine $\mu$ uniquely.
They therefore determine $\Phi$.
\end{theorem}
\end{samepage}

\begin{corollary}[Reconstruction from the endpoint coefficients]
\label{cor:endpoint-reconstruction}
At all scales $r$ and levels $N$, restrict to the tests $f^{r,N}_{x,b}$
with $\langle x,b\rangle\not\equiv0\pmod N$.  For $k\ge4$, the values
of $\mu_k^{(1)}$ and $\mu_k^{(k-1)}$ on these tests determine the
weight-$k$ component of $\Phi$ uniquely among symbols with the action
\eqref{eq:old-action}.  For $k=2,3$, the values of $\mu_k^{(1)}$ alone
suffice.
\end{corollary}

A nonzero pairing $\langle x,b\rangle$ modulo $N$ allows a rational
change of box into $\G(\Z_p)$ for some $p\mid N$, where the local
reciprocity law applies.  For $j\ne2$ one can divide out the finite-level
smoothing scalar.  For uniqueness, consider the difference of two
distributions with the same endpoint data.  Adding the second column
of a box to the first acts on
$f^{r,N}_{x,b}$ by a character.  The nonzero pairings give exactly the
nontrivial eigenspaces of this translation.  Vanishing there forces
the tested coefficients to be invariant under the translation and,
by diagonal conjugation, under every rational lower unipotent matrix.
In weight $k\ge4$ the first Manin relation exchanges the two endpoint
coefficients, giving invariance under the upper unipotents as well.
Their torus characters then force both to vanish.  The cyclic relation
gives a linear dependence among the interior coefficients, which the
diagonal torus separates by their distinct characters
(Theorem~\ref{thm:Fourier-Manin-uniqueness}).

At level one, Hirose, Sato, and Tasaka obtain a complete family of
linear relations among $\mathbb E_k$ and the products
$\mathbb E_{2r}\mathbb E_{k-2r}$ by partial fraction decomposition
\cite[Theorem~1 and Corollary~2]{HiroseSatoTasaka}.  Here the relations
arise as coefficients of a single cyclic Manin identity evaluated on
the integral lattice.  In even weight $k\ge8$, the combinations of
these coefficients that cancel $\mathbb E_k$ give all relations among
the unordered products.  Their space has
dimension $\lfloor(k-8)/6\rfloor$
(Theorem~\ref{thm:journal-level-one-relations}).

At prime level, Borisov, Gunnells, and Popescu used an auxiliary
elliptic curve to obtain equations for $X_1(p)$ in weight-one
Eisenstein coordinates \cite[Theorem~4.4 and Corollary~4.5]{BGP}.
For the full-level curve $X(p)$, Khuri-Makdisi observes that the
weight-one Eisenstein embedding admits defining equations of degrees
two and three \cite[end of \S5]{KhuriMakdisi}.  We prove that, for
every odd prime $p$, the quadrics obtained from the cyclic Manin identity
by eliminating the holomorphic weight-two coordinates cut out this
embedding scheme-theoretically.  Consequently every quadratic relation modulo
weight-two Eisenstein series follows from the cyclic equations after
saturation by the weight-one variables
(Theorem~\ref{thm:prime-weight-two-saturation}).  For $X(p)$ the equations on the separate
cyclic subgroups must recover one elliptic curve carrying the full
$p$-torsion with its prescribed Weil pairing.  Appendix~\ref{app:prime-weight-two}
proves this for $p\ge5$ and computes the tangent space at the cusps.
At $p=3$ the weight-one embedding identifies $X(3)$ with the whole
projective line.  The corresponding saturation assertion in arbitrary
level $N\ge3$ and weight is Conjecture~\ref{conj:R-relation-ideal}.

Khuri-Makdisi describes torsion Eisenstein series through Laurent
expansions of elliptic functions \cite[\S2]{KhuriMakdisi}.
Khuri-Makdisi and Raji derive the product relations modulo Eisenstein
series from Petersson pairings and periods
\cite[Theorem~2.8]{KhuriMakdisiRaji}.  For series defined by additive
characters, Zhang proves the weight-two identity
\cite[Corollary~1.4.9]{Zhang}, and Brunault gives the explicit
Eisenstein term in every weight for nonzero parameters
\cite[Theorem~1]{BrunaultBG}.  Brunault identifies Theorem~3.1 of
the 2006 version of this article \cite{Pasol2006} as the analogue for
translated series.  He also asks about level one and the zero
parameter \cite[\S6]{BrunaultBG}.  In the translated family used here,
omitting the singular lattice term gives the level-one identity
\eqref{eq:journal-full-triangle}.

There are also geometric and motivic constructions of related symbols.
Bergeron, Charollois, and Garcia obtain modular forms from homology by
theta lifting \cite{BCG}; Xu gives an explicit formula in rank two
\cite{Xu}.  In weight two, Branchereau identifies products of weight-one
Eisenstein series with periods of such a lift
\cite[Theorem~1.4]{Branchereau}.  Sharifi and Venkatesh construct motivic
cocycles whose torsion specializations give symbols of Siegel units,
and prove their Hecke equivariance \cite{SharifiVenkatesh}.

Sections~\ref{sec:coefficients}--\ref{sec:symbol} give the analytic
construction.  Section~\ref{sec:Colmez} proves the finite Fourier
identity and uniqueness; Section~\ref{sec:BK} gives the reconstruction.
Section~\ref{sec:level-one} treats relations among Eisenstein products.
Appendix~\ref{app:prime-weight-two} proves the prime-level theorem.

\subsection*{Acknowledgements}
The author is grateful to Glenn Stevens for his guidance and for the ideas
that led to the original construction.

\section{Test functions and Eisenstein--Kronecker series}\label{sec:coefficients}

\subsection{Test functions and distributions}\label{sec:test-functions}

We use the convention of \cite[\S2.1]{BPPS}.  Let $V$ be a
finite-dimensional $\Q$-vector space.  A \emph{lattice} in $V$ is a
finitely generated additive subgroup that spans $V$ over $\Q$.  A function
$\varphi:V\to\Z$ has bounded support if its support is contained in a
lattice; it is locally constant for the lattice topology if its group of
periods contains a lattice.  Define the BPPS test-function space
\begin{equation}\label{eq:BPPS-test-space}
 \cS(V)=\{\varphi:V\to\Z:
       \varphi\text{ has bounded support and is locally constant}\}.
\end{equation}
Its elements are called \emph{test functions}.
If $A$ is an abelian group, put
\[
 \Dist(V,A)=\operatorname{Hom}_{\Z}(\cS(V),A).
\]
For $U\subset V$ with characteristic function $[U]\in\cS(V)$, write
$\nu(U)=\nu([U])$.  The characteristic functions of lattice cosets generate
$\cS(V)$ over $\Z$.

The canonical identification
\begin{equation}\label{eq:tensor-test}
 \cS(\Q^2)\otimes_{\Z}\cS(\Q^2)
 \xrightarrow{\sim}
 \cS(\Mat(\Q))
\end{equation}
sends $\varphi_1\otimes\varphi_2$ to the function whose value on a matrix with columns $v_1,v_2$ is $\varphi_1(v_1)\varphi_2(v_2)$.

\subsection{Group-action conventions}

Matrices act on $\mathbf P^1(\Q)$ by fractional linear transformations;
this induces their action on $\Delta_0$.  Fix $\gamma\in\G(\Q)$ and a test
function $\varphi$ on $\Mat(\Q)$.  Put
\[
 (\gamma_*\varphi)(A)=\varphi(A\gamma).
\]
If $W$ is a left $\G(\Q)$-module, the induced action on $\Dist(\Mat(\Q),W)$ is fixed by
\begin{equation}\label{eq:dist-action}
 (\gamma\cdot\nu)(\varphi)
 =\gamma\cdot_W\nu((\gamma^{-1})_*\varphi).
\end{equation}

For the formal modular coefficient modules below, right multiplication on the matrix
variable mixes the two columns.  We normalize the induced action on a
distribution so that on characteristic functions $[U]\in\cS(\Mat(\Q))$
\begin{equation}\label{eq:old-action}
 (\gamma\cdot\nu)(U)(X,Y)
 =\det(\gamma)\,
 \nu(U\gamma)\bigl((X,Y)\gamma\bigr).
\end{equation}
The determinant and substitution in \eqref{eq:old-action} are part of the
coefficient action.  Formula~\eqref{eq:old-action} is algebraic: no slash
operator on the variable $\tau$ is involved, so it defines the action for
either sign of the determinant.  Section~\ref{sec:torus} verifies directly
the covariance needed for the construction.

\subsection{Modular and nearly holomorphic forms}

Let $\Gamma\subset\SL(\Z)$ be a congruence subgroup.  We write
\[
 M_k(\Gamma),\quad S_k(\Gamma),\quad \cE_k(\Gamma)
\]
for the spaces of holomorphic modular forms, cusp forms, and Eisenstein forms of weight $k$.  We use the Hecke-stable decomposition
\begin{equation}\label{eq:Eis-decomp}
 M_k(\Gamma)=S_k(\Gamma)\oplus\cE_k(\Gamma)
\end{equation}
whenever it is available in the stated weight and level.  Independently of a splitting, put
\[
 \cC_k(\Gamma)=M_k(\Gamma)/\cE_k(\Gamma).
\]

Put $y=\operatorname{Im}(\tau)$.  A nearly holomorphic modular form of
weight $k$ and depth at most $r$ is a function of the form
\[
 f(\tau)=\sum_{j=0}^r f_j(\tau)y^{-j},
\]
where the $f_j$ are holomorphic on $\Hh$.  For $f\ne0$, its depth is the
degree of this polynomial in $y^{-1}$.
The function $f$ transforms with
weight $k$ under $\Gamma$ and has polynomial growth at every cusp.  Write
\[
 \delta_k=\frac{\partial}{\partial\tau}+\frac{k}{\tau-\bar\tau},
 \qquad
 \varepsilon=(\tau-\bar\tau)^2\frac{\partial}{\partial\bar\tau}.
\]
Equivalently, with
$\partial_{\bar\tau}$ interpreted by treating $\tau$ and $\bar\tau$ as
independent variables, it satisfies
$\varepsilon^{r+1}f=0$.  We denote this space by $N_k^r(\Gamma)$.  This is
Shimura's notion of a nearly holomorphic modular form
\cite[\S2]{ShimuraNH}; see also \cite[Chapter~10]{Hida}.
Put
\[
 E_2^*(\tau)=1-24\sum_{n\geq1}\sigma_1(n)e^{2\pi in\tau}
               -\frac{3}{\pi y},
\]
the level-one nearly holomorphic Eisenstein series of weight two.  For
$j\geq1$, write
\[
 \delta_w^j=\delta_{w+2j-2}\circ\cdots\circ\delta_{w+2}\circ\delta_w,
 \qquad \delta_w^0=\operatorname{id}.
\]

For $k\ge1$ and $r\ge0$, put
\[
 J(k,r)=\{j\in\Z:0\le j\le r,\ k-2j\ge1\}.
\]
\[
 X_k^r=
 \begin{cases}
  \C\,\delta_2^{(k-2)/2}E_2^*,
     &k\text{ even and }r\ge k/2,\\
  0,&\text{otherwise}.
 \end{cases}
\]

\begin{theorem}[Nearly holomorphic decomposition]
\label{thm:H-package}
For $k\ge1$ and $r\ge0$,
\begin{equation}\label{eq:Shimura-decomp}
 N_k^r(\Gamma)=
 \bigoplus_{j\in J(k,r)}
 \delta_{k-2j}^{j}M_{k-2j}(\Gamma)\oplus X_k^r.
\end{equation}
Moreover, $N_0^r(\Gamma)=\C$.
\end{theorem}

\begin{proof}
This is \cite[Theorem~4.2]{PitaleSahaSchmidt}, whose raising operator is
$R_k=2i\delta_k$.  For even $k$, their exceptional condition
$2\le k<2+2r$ is equivalent to $r\ge k/2$.
\end{proof}

Define the algebraic holomorphic projection
\[
 H_{k,\Gamma}^{\mathrm{alg}}:
 N_k^r(\Gamma)\longrightarrow M_k(\Gamma)
\]
as projection onto the $j=0$ summand in
\eqref{eq:Shimura-decomp}.  It is independent of $r$.
Suppressing $\Gamma$ when unambiguous, one has
\[
 H_k^{\mathrm{alg}}|_{M_k(\Gamma)}=\operatorname{id},
 \qquad
 H_{k+2}^{\mathrm{alg}}\circ\delta_k=0.
\]
Projection commutes with restriction to a congruence subgroup and,
for $\gamma\in\SL(\Z)$, satisfies
\[
 H_{k,\gamma^{-1}\Gamma\gamma}^{\mathrm{alg}}(f|_k\gamma)
 =H_{k,\Gamma}^{\mathrm{alg}}(f)|_k\gamma.
\]
Indeed, restriction preserves the decomposition, while
\[
 \delta_w(f|_w\gamma)=(\delta_wf)|_{w+2}\gamma,
 \qquad E_2^*|_2\gamma=E_2^*.
\]
The assertions follow from uniqueness in
\eqref{eq:Shimura-decomp}.  Composing with
$M_k(\Gamma)\to\cC_k(\Gamma)$ requires no choice of an
Eisenstein complement.

We shall apply this projection to products $E_{m,a}E_{n,b}$, $m,n\ge1$,
of the Eisenstein--Kronecker series defined in Section~\ref{sec:EK},
including their regularized values at zero.  Each factor is holomorphic
unless its weight is two, when it has depth one.  Thus each product has
weight $k=m+n\ge2$ and is a polynomial of degree at most two in $y^{-1}$.
In this range, $X_k^r$ can be nonzero only for $k=2$ or $k=4$,
when it is respectively $\C E_2^*$ or $\C\delta_2E_2^*$.
The slash compatibility above concerns $\SL(\Z)$; the rational
parameter action is defined separately by \eqref{eq:old-action}.

\subsection{Power series with modular-form coefficients}
\label{sec:power-series-modules}

For $k<0$, set
\[
 M_k(\Gamma)=S_k(\Gamma)=\cE_k(\Gamma)=\cC_k(\Gamma)=0.
\]
For $m\in\Z$, consider formal series whose coefficient of $X^iY^j$
has modular weight $m+i+j$.  Write
\begin{equation}\label{eq:formal-modular-module}
 \scrM_m(\Gamma)=
 \prod_{(i,j)\in\Z_{\geq0}^2}
 M_{m+i+j}(\Gamma)X^iY^j.
\end{equation}
No convergence in $X,Y$ is required.  Homogeneous degree $n$ has
modular weight $m+n$.

We use the following notation:
\begin{align*}
 \scrN_m^r(\Gamma)
   &=\prod_{i,j\geq0}N_{m+i+j}^r(\Gamma)X^iY^j,\\
 \scrE_m(\Gamma)
   &=\prod_{i,j\geq0}\cE_{m+i+j}(\Gamma)X^iY^j,\\
 \scrM_m^{\mathrm{cusp}}(\Gamma)
   &=\prod_{i,j\geq0}S_{m+i+j}(\Gamma)X^iY^j,\\
 \scrC_m(\Gamma)
   &=\prod_{i,j\geq0}\cC_{m+i+j}(\Gamma)X^iY^j.
\end{align*}
We call their elements, respectively, nearly holomorphic coefficient series
of depth at most $r$, Eisenstein-coefficient series, cusp-form-coefficient
series, and cuspidal quotient coefficient series.
Coefficientwise quotienting gives the canonical identification
\[
 \scrC_m(\Gamma)=\scrM_m(\Gamma)/\scrE_m(\Gamma).
\]
After a choice of the decompositions \eqref{eq:Eis-decomp}, this quotient
is represented by $\scrM_m^{\mathrm{cusp}}(\Gamma)$; none of the canonical
statements below depends on that choice.

Passing to the filtered colimit over congruence subgroups gives
\[
 \scrM_m,\quad\scrN_m^r,\quad\scrE_m,\quad
 \scrM_m^{\mathrm{cusp}},\quad\scrC_m.
\]
For $F=\sum_{i,j}f_{i,j}X^iY^j\in\scrN_m^r$, define
\begin{equation}\label{eq:coefficientwise-H}
 \begin{aligned}
 \widetilde{\cH}(F)
 &=\sum_{i,j\geq0}
 H_{m+i+j}^{\mathrm{alg}}(f_{i,j})X^iY^j
 \in\scrM_m,\\
 \cH(F)
 &=\sum_{i,j\geq0}
 \pi_{\cC_{m+i+j}}H_{m+i+j}^{\mathrm{alg}}(f_{i,j})X^iY^j
 \in\scrC_m.
 \end{aligned}
\end{equation}
where $\pi_{\cC_k}:M_k\to\cC_k$ is the quotient map.  We call $\cH$ the
\emph{coefficientwise cuspidal holomorphic projection};
$\widetilde{\cH}$ is the coefficientwise holomorphic projection before
passing to the Eisenstein quotient, and $\cH(F)$ is its image in $\scrC_m$.
For any bivariate
coefficient series $G(X,Y)$, write $G^{[n]}$ for its homogeneous
total-degree-$n$ part.

\subsection{Eisenstein--Kronecker series}\label{sec:EK}

We recall the analytic continuation and torsion specialization of the
Kronecker series.  The resulting series $E_V$ are the one-variable factors
used in the definition of $\mu$ and in the Fourier comparison.

Fix $\tau\in\Hh$ and put $L_\tau=\Z\tau+\Z$.  If $k\geq1$, define initially in a right half-plane
\begin{equation}\label{eq:phi-def}
 \phi_k(z,\tau;s)
 =\sum_{\omega\in L_\tau}^{\prime}
 (z+\omega)^{-k}|z+\omega|^{-s},
\end{equation}
where the term $\omega=-z$ is omitted if $z\in L_\tau$.

\begin{theorem}[Analytic continuation and differentiation]\label{thm:analytic-cont}
The series \eqref{eq:phi-def} converges normally for $\operatorname{Re}(s)>2-k$, has meromorphic continuation to the $s$-plane, and is holomorphic at $s=0$.  With
\[
 \phi_k(z,\tau)=\phi_k(z,\tau;0),
\]
one has, away from $L_\tau$,
\[
 \partial_z\phi_k(z,\tau)=-k\phi_{k+1}(z,\tau),
 \qquad
 \phi_k(-z,\tau)=(-1)^k\phi_k(z,\tau).
\]
The continuation, differentiation, and parity statements are locally uniform in $(z,\tau)$ on compact subsets avoiding the polar locus.
\end{theorem}

\begin{proof}
This is the classical continuation of the Kronecker double series;
see \cite[Chapter~VIII, \S13]{Weil} and \cite[\S1.3.1]{Colmez}.
In the notation of \cite[\S1.1]{BannaiKobayashi}, our series is
$K_k^*(z,0,k+s/2;L_\tau)$.  The Mellin integral after Poisson summation
expresses it, up to a reciprocal gamma factor, as a sum of two integrals
of exponentially decreasing theta series.  Since $k>0$, there are no
polar terms.  These integrals and their derivatives converge uniformly
on compact sets avoiding the polar locus.  Normal convergence in
$\operatorname{Re}(s)>2-k$ permits termwise differentiation with respect to
$z$, with $\bar z$ held fixed.
For general $s$ one has
\[
 \partial_z\bigl((z+\omega)^{-k}|z+\omega|^{-s}\bigr)
 =-(k+s/2)(z+\omega)^{-k-1}|z+\omega|^{-s};
\]
evaluating the continued identity at $s=0$ gives the displayed derivative
formula.  Reindexing $\omega\mapsto-\omega$ gives parity.  The same
integral estimates justify differentiation of the continued identity.
\end{proof}

For $a=(a_1,a_2)\in\Q^2$, put $\tau_a=a_1\tau+a_2$ and
\[
 E_{k,a}(\tau)=\phi_k(\tau_a,\tau).
\]
The \emph{exact level} of $a\notin\Z^2$ is the order of its class in
$(\Q/\Z)^2$; equivalently, it is the least integer $N\geq1$ such that
$Na\in\Z^2$.

\begin{theorem}[Torsion specialization]\label{thm:torsion-specialization}
Let $a\in\Q^2\setminus\Z^2$ and $k\geq1$.
\begin{enumerate}[label=\textup{(\roman*)}]
 \item For $\gamma\in\SL(\Z)$,
 \[
  E_{k,a}(\gamma\tau)j(\gamma,\tau)^{-k}=E_{k,\gamma^ta}(\tau).
 \]
 \item If $k\neq2$, then $E_{k,a}$ is holomorphic in $\tau$.
 \item If $k=2$, then $E_{2,a}-E_{2,0}$ is holomorphic and equals $\wp(\tau_a,\tau)$.
 \item If $N$ is the exact level of $a$, then $E_{k,a}$ is Eisenstein on
 $\Gamma(N)$ for $k\ne2$, while $E_{2,a}-E_{2,0}$ belongs to the holomorphic
 Eisenstein subspace of weight two on $\Gamma(N)$.
\end{enumerate}
\end{theorem}

\begin{proof}
If $\gamma=\begin{psmallmatrix}p&q\\r&s\end{psmallmatrix}$, then
$L_{\gamma\tau}=j(\gamma,\tau)^{-1}L_\tau$ and
\[
 a_1\gamma\tau+a_2
 =j(\gamma,\tau)^{-1}\tau_{\gamma^ta}.
\]
Changing variables in the continued double series proves (i).  Assertions
(ii)--(iv) follow from the torsion specialization of Kronecker series;
in our normalization these are the series $G_k$ of
\cite[Definition~2.1 and Proposition~2.10]{KhuriMakdisi}.
More explicitly, the theta formula
\eqref{eq:phi1-theta} below shows holomorphy in weight one, differentiation
gives the weights $k\geq3$, and in weight two the sole nonholomorphic summand
is the scalar series $E_{2,0}$; subtraction leaves the Weierstrass function.
Proposition~2.10 of loc.\ cit. identifies their span with the Eisenstein
subspace, using differences from $E_{2,0}$ in weight two.  At level one
the assertion is the classical lattice definition of Eisenstein series.
\end{proof}

\subsection{Eisenstein--Kronecker series on test functions}

For a subset $V\subset\Q^2$ with $[V]\in\cS(\Q^2)$, define
\[
 V_\tau=\{a_1\tau+a_2:(a_1,a_2)\in V\}
\]
and
\[
 \phi_{k,V}(z,\tau;s)
 =\sum_{\omega\in V_\tau}^{\prime}
 (z+\omega)^{-k}|z+\omega|^{-s}.
\]

\begin{proposition}\label{prop:compact-cont}
The function $\phi_{k,V}(z,\tau;s)$ has meromorphic continuation in $s$ and is regular at $s=0$.  The continuation is independent of a presentation of $V$ as a disjoint union of rational lattice cosets.  Moreover, the assignment
\[
 [V]\longmapsto \phi_{k,V}(z,\tau;0)
\]
extends linearly to a distribution in the $V$-variable.
\end{proposition}

\begin{proof}
Write $V$ as a finite disjoint union of cosets $a+r\Z^2$, with one common
$r\in\Q^\times$.  For one coset, direct rescaling in the half-plane of
absolute convergence gives
\[
 \phi_{k,a+r\Z^2}(z,\tau;s)
 =r^{-k}|r|^{-s}
 \phi_k\left(\frac{z+\tau_a}{r},\tau;s\right).
\]
Theorem~\ref{thm:analytic-cont} therefore continues each summand and proves
regularity at zero.  Any two decompositions have a common lattice refinement.
Finite additivity holds before continuation, hence after continuation by the
identity theorem.  This proves independence and gives the asserted linear map
on characteristic functions, hence on all test functions.
\end{proof}

\section{Quadratic identities and holomorphic projection}\label{sec:quadratic}

\subsection{Theta functions and the fundamental three-term identity}\label{sec:theta}

The theta-function identity below supplies the second Manin relation.  Its
regularized forms are used in Sections~\ref{sec:generating} and
\ref{sec:Manin}.

Put $q=e^{2\pi i\tau}$ and use the odd Jacobi theta function
\[
 \theta(z,\tau)=\frac1i\sum_{n\in\Z}(-1)^n
 q^{(2n+1)^2/8}e^{(2n+1)\pi iz}.
\]
Let $\sigma(z,\tau)$, $\zeta(z,\tau)$, and $\wp(z,\tau)$ be the
Weierstrass sigma, zeta, and $\wp$ functions for $L_\tau$, normalized by
$\zeta=\partial_z\log\sigma$ and $\wp=-\partial_z\zeta$.  These
normalizations and the formulas below are those of
\cite[Chapters~IV--V]{Chandrasekharan}.

\begin{lemma}[Theta and Weierstrass addition formulas]\label{lem:theta-addition}
There is a quasi-period $\eta_1(\tau)$ such that
\[
 \sigma(z,\tau)=e^{\eta_1z^2}\frac{\theta(z,\tau)}{\theta'(0,\tau)},
\]
and, meromorphically in $z_1,z_2$,
\begin{align*}
 \wp(z_1)-\wp(z_2)
 &=-\frac{\theta(z_1+z_2)\theta(z_1-z_2)\theta'(0)^2}
 {\theta(z_1)^2\theta(z_2)^2},\\
 \frac{\wp'(z_1)-\wp'(z_2)}{\wp(z_1)-\wp(z_2)}
 &=-2\bigl(\partial_z\log\theta(z_1)+\partial_z\log\theta(z_2)
 -\partial_z\log\theta(z_1+z_2)\bigr).
\end{align*}
\end{lemma}

\begin{proof}
The first two formulas are
\cite[Chapter~V, \S2, Theorem~2]{Chandrasekharan}.  Logarithmically
differentiate the second formula with respect to $z_1$ and then with respect
to $z_2$.  Add the resulting identities and use that $\theta$ is odd, so
$\theta'/\theta$ is odd.  This gives the final formula, first off the
divisors and then everywhere by meromorphic continuation.
\end{proof}

\begin{theorem}[Frobenius--Stickelberger identity]\label{thm:FS}
If $z_1+z_2+z_3=0$ and no $z_i$ lies in $L_\tau$, then
\begin{equation}\label{eq:FS}
 \bigl(\phi_1(z_1,\tau)+\phi_1(z_2,\tau)+\phi_1(z_3,\tau)\bigr)^2
 =\wp(z_1,\tau)+\wp(z_2,\tau)+\wp(z_3,\tau).
\end{equation}
After imposing \(z_3=-z_1-z_2\), both sides are meromorphic functions of
\((z_1,z_2)\) and agree on the complement of their polar divisors; hence they
agree identically as meromorphic functions.
\end{theorem}

\begin{proof}
The Kronecker theta formula is
\begin{equation}\label{eq:phi1-theta}
 \phi_1(z,\tau)=\partial_z\log\theta(z,\tau)
 +2\pi i\frac{z-\bar z}{\tau-\bar\tau}
\end{equation}
by \cite[Chapter~VIII, \S13]{Weil}; it also appears in
\cite[Chapter~II]{deShalit}.  The affine correction in
\eqref{eq:phi1-theta} cancels when $z_1+z_2+z_3=0$.  On the other hand, the
Weierstrass addition law
\[
 \wp(z_1+z_2)=\frac14
 \left(\frac{\wp'(z_1)-\wp'(z_2)}{\wp(z_1)-\wp(z_2)}\right)^2
 -\wp(z_1)-\wp(z_2)
\]
together with Lemma~\ref{lem:theta-addition} gives
\[
 \wp(z_1)+\wp(z_2)+\wp(z_1+z_2)
 =\bigl(h(z_1)+h(z_2)-h(z_1+z_2)\bigr)^2,
\]
where $h=\partial_z\log\theta$.  Since $\wp$ is even and $h$ is odd,
putting $z_3=-z_1-z_2$ proves \eqref{eq:FS}.  The meromorphic extension
removes the temporary restrictions used in the divisions above.
\end{proof}

For a $C^\infty$ function $f(z,\bar z,\tau,\bar\tau)$, put
\[
 \Tay_{z_0}f(\tau,X)
 =\sum_{n\geq0}\frac{\partial_z^nf(z_0,\tau)}{n!}X^n.
\]
If $f$ is meromorphic in $z$ and regular at $z_0$, this is its ordinary
Taylor expansion.  At a pole it is instead a Laurent germ and is not covered
by the notation above.  For $\phi_1$ at a regular point it is the formal
Taylor series for $\partial_z$, with $\bar z$ held fixed.

\begin{theorem}[Formal fundamental identity]\label{thm:fundamental-identity}
Suppose $z_1+z_2+z_3=0$ and $z_i\notin L_\tau$ for
$i=1,2,3$.  Then, as an identity of analytic germs and hence of formal
power series,
\begin{align}
 &\bigl(
 \Tay_{z_1}\phi_1(\tau,X)
 +\Tay_{z_2}\phi_1(\tau,Y)
 +\Tay_{z_3}\phi_1(\tau,-X-Y)
 \bigr)^2 \notag\\
 &\qquad=
 \Tay_{z_1}\wp(\tau,X)
 +\Tay_{z_2}\wp(\tau,Y)
 +\Tay_{z_3}\wp(\tau,-X-Y).
 \label{eq:fundamental-identity}
\end{align}
\end{theorem}

\begin{proof}
Write $h=\partial_z\log\theta$ and
$A(z)=(z-\bar z)/(\tau-\bar\tau)$.  Formula
\eqref{eq:phi1-theta} gives
\begin{align*}
 &\Tay_{z_1}\phi_1(X)+\Tay_{z_2}\phi_1(Y)
 +\Tay_{z_3}\phi_1(-X-Y)\\
 &\qquad=\Tay_{z_1}h(X)+\Tay_{z_2}h(Y)+\Tay_{z_3}h(-X-Y),
\end{align*}
because both the constant and linear terms contributed by $A$ cancel.
The right-hand side is a sum of genuine meromorphic Taylor series.  Substitute
$z_1+x,z_2+y,z_3-x-y$ in the theta form of the
Frobenius--Stickelberger identity proved above and compare Taylor
coefficients.  This proves \eqref{eq:fundamental-identity}.  The hypothesis
$z_i\notin L_\tau$ makes each translated expression an analytic germ.  The
argument uses formal translation for the derivation $\partial_z$ and does
not require $\phi_1$ to be holomorphic in $z$.
\end{proof}

\begin{remark}[Lattice parameters]\label{rem:fundamental-lattice-points}
At a lattice point the identity is an identity of Laurent germs.
Subtracting the polar terms gives the power-series identities of
Proposition~\ref{prop:regularized-fundamental}.
\end{remark}

\subsection{Generating series and regularization}\label{sec:generating}

We define the series $E_a$, $E_V$, and the regularized series $F$ attached
to the zero lattice coset.  They enter the definition of $\mu$ and the proof
of the Manin relations when a torsion parameter is zero.

For $a\in\Q^2\setminus\Z^2$, define the generating series
$E_a(\tau,X)$ and its coefficients $E_{n+1,a}(\tau)$ by
\begin{equation}\label{eq:Ea}
 E_a(\tau,X):=\Tay_{\tau_a}\phi_1(\tau,X)
 =\sum_{n\geq0}(-1)^nE_{n+1,a}(\tau)X^n.
\end{equation}
For a subset $V\subset\Q^2$ with $[V]\in\cS(\Q^2)$, define
\begin{equation}\label{eq:EV}
 E_V(\tau,X)
 :=\sum_{n\geq0}(-1)^n\phi_{n+1,V}(0,\tau)X^n,
\end{equation}
where the zero summand is omitted in every coefficient.  For
$a\notin\Z^2$, the choice $V=a+\Z^2$ recovers \eqref{eq:Ea}:
$E_{a+\Z^2}(\tau,X)=E_a(\tau,X)$.  If $0\notin V$,
this is the formal Taylor series $\Tay_0\phi_{1,V}$ for $\partial_z$.
If $0\in V$, it is
the Taylor series of the regular part obtained by subtracting the polar term
of the unique refined lattice coset containing zero; this is made explicit in
Corollary~\ref{cor:zero-coset}.

For a one-variable formal series $G(X)=\sum_{n\geq0}g_nX^n$, define its
formal primitive with zero constant term by
\begin{equation}\label{eq:formal-integral}
 \Int_XG=\sum_{n\geq0}g_n\frac{X^{n+1}}{n+1}.
\end{equation}
If $G=\sum_{n\geq0}g_nX^n$ has weight shift $m$, so that $g_n$ has
weight $m+n$, define the coefficientwise raising operator by
\begin{equation}\label{eq:series-raising}
 \delta G=\sum_{n\geq0}\delta_{m+n}(g_n)X^n.
\end{equation}

\begin{proposition}[Scaling and parity]\label{prop:scaling}
For $\alpha\in\Q^\times$,
\[
 E_{\alpha V}(\tau,\alpha X)=\alpha^{-1}E_V(\tau,X),
 \qquad
 E_{-V}(\tau,-X)=-E_V(\tau,X).
\]
The same identities hold for arbitrary test functions by linearity.
\end{proposition}

\begin{proof}
For $\operatorname{Re}(s)$ large, substitution $\omega=\alpha\omega'$ gives
\[
 \phi_{n,\alpha V}(0,\tau;s)
 =\alpha^{-n}|\alpha|^{-s}\phi_{n,V}(0,\tau;s).
\]
Continue to $s=0$, multiply the coefficient of $X^{n-1}$ by
$(\alpha X)^{n-1}$, and sum; this gives the factor $\alpha^{-1}$.
Parity follows in the same way from $\omega\mapsto-\omega$.  Linearity gives
the assertion for test functions.
\end{proof}

\subsubsection{The integral lattice}

The series at the zero torsion point is regularized by subtracting its polar part:
\begin{equation}\label{eq:regularized-F}
 F(\tau,X)=E_{\Z^2}(\tau,X)
 =\left.\Tay_z\left(\phi_1(z,\tau)-\frac1z\right)(X)\right|_{z=0}.
\end{equation}

\begin{proposition}[Regularized limit identities]\label{prop:regularized}
The series $F$ is odd in $X$, and its coefficient of $X^n$ belongs to
$N_{n+1}^1(\SL(\Z))$.  With the zero summand omitted, it satisfies the
same scaling and lattice-refinement identities as the series $E_V$.
Its square admits the identity
\begin{equation}\label{eq:F-square}
 F(X)^2
 =-\partial_XF(X)-2\frac{F(X)}{X}
 +4\pi i\,\delta\Int_XF(X)-3E_{2,0}.
\end{equation}
\end{proposition}

\begin{proof}
By \eqref{eq:phi1-theta}, the difference $\phi_1(z,\tau)-z^{-1}$ is
real analytic at zero.  The derivative identities of
Theorem~\ref{thm:analytic-cont} identify its Taylor coefficients with
those of $E_{\Z^2}$.  Pairing $\omega$ with $-\omega$ proves oddness.
Theorem~\ref{thm:torsion-specialization} gives modularity and growth;
only $E_{2,0}$ is nonholomorphic.  Scaling and refinement follow from
the corresponding lattice sums, first in their domain of convergence
and then by continuation to $s=0$.

For the square identity, take the regular part of the logarithmic theta
derivative:
\[
 P(X)=\frac{\theta'(X)}{\theta(X)}-\frac1X,
 \qquad F(X)=P(X)+\frac{2\pi iX}{\tau-\bar\tau}.
\]
The heat equation $\theta''=4\pi i\partial_\tau\theta$ gives
\[
 P^2+P'+2P/X
 =4\pi i\partial_\tau\log\theta'(0)
   +4\pi i\partial_\tau\Int_XP.
\]
Comparison of the linear coefficients gives
\[
 4\pi i\partial_\tau\log\theta'(0)
 =\frac{\theta'''(0)}{\theta'(0)}
 =-3E_{2,0}-\frac{6\pi i}{\tau-\bar\tau}.
\]
Substitute $P=F-2\pi iX/(\tau-\bar\tau)$.  The terms in
$(\tau-\bar\tau)^{-2}X^2$ cancel, leaving
\[
 F^2+F'+2F/X
 =-3E_{2,0}+4\pi i\left(\partial_\tau\Int_XF
                    +\frac{XF}{\tau-\bar\tau}\right).
\]
The expression in parentheses is $\delta\Int_XF$, since the coefficient
of $X^n$ in $\Int_XF$ has weight $n$.  This proves \eqref{eq:F-square}.
\end{proof}

\begin{corollary}[The zero lattice coset]\label{cor:zero-coset}
For $r\in\Q^\times$, the series attached to the zero coset modulo
$r\Z^2$ is
\begin{equation}\label{eq:scaled-zero-coset}
 E_{r\Z^2}(X)=r^{-1}F(X/r).
\end{equation}
If $[V]\in\cS(\Q^2)$ and $V$ is written as a disjoint union of rational lattice cosets and one of them is the zero coset, that summand is interpreted by
\eqref{eq:scaled-zero-coset}.  The resulting series is independent of the
chosen common refinement.
\end{corollary}

\begin{proof}
Apply Proposition~\ref{prop:scaling} to $V=\Z^2$ and replace its formal
variable by $X/r$.  Independence follows from
Proposition~\ref{prop:compact-cont}: two decompositions have a common
refinement, and the same scaling identity applies to the unique refined
coset containing zero.
\end{proof}

For $a\in\Q^2\setminus\Z^2$, put
\[
 W_a(T)=\Tay_{\tau_a}\wp(\tau,T),
 \qquad
 W_0(T)=\Tay_0\bigl(\wp(z,\tau)-z^{-2}\bigr)(T).
\]

\begin{proposition}[Regularized fundamental identities]
\label{prop:regularized-fundamental}
The fundamental identity has the following power-series forms at its polar
parameters.
\begin{enumerate}[label=\textup{(\roman*)}]
 \item If $a\in\Q^2\setminus\Z^2$ and
 \[
  P_a(X,Y)=E_a(X)+F(Y)+E_{-a}(-X-Y),
 \]
 then $P_a$ is divisible by $Y$ and
 \begin{equation}\label{eq:one-integral}
  P_a^2=W_a(X)+W_0(Y)+W_{-a}(-X-Y)-2P_a/Y.
 \end{equation}
 This is the regularized identity when exactly one of the three base points
 lies in $L_\tau$; the other placements follow by cyclic permutation.
 \item Put
 \[
  P_0(X,Y)=F(X)+F(Y)+F(-X-Y),
  \qquad
  Q(X,Y)=\frac1X+\frac1Y-\frac1{X+Y}.
 \]
 Then $P_0$ is divisible by $XY(X+Y)$, so $P_0Q$ is a power series, and
 \begin{equation}\label{eq:all-integral-square}
  P_0^2=W_0(X)+W_0(Y)+W_0(-X-Y)-2P_0Q.
 \end{equation}
 Equivalently,
 \begin{align}\label{eq:all-integral}
  &F(X)F(Y)+F(Y)F(-X-Y)+F(-X-Y)F(X)\\
  &\quad=\frac12\bigl(W_0(X)+W_0(Y)+W_0(-X-Y)\bigr)\notag\\
  &\qquad\quad-\frac12\bigl(F(X)^2+F(Y)^2+F(-X-Y)^2\bigr)-P_0Q.\notag
 \end{align}
 This is the regularized identity when all three base points lie in
 $L_\tau$.
\end{enumerate}
\end{proposition}

\begin{proof}
For (i), apply the meromorphic Frobenius--Stickelberger identity to
\[
 z_1=\tau_a+X,\qquad z_2=Y,\qquad
 z_3=-\tau_a-X-Y.
\]
At the polar point use the Laurent decompositions
\[
 \phi_1(Y)=Y^{-1}+F(Y),\qquad
 \wp(Y)=Y^{-2}+W_0(Y).
\]
After the $Y^{-2}$ terms cancel, the remaining identity is
\eqref{eq:one-integral}.  Parity gives
$P_a(X,0)=E_a(X)+E_{-a}(-X)=0$, hence $Y\mid P_a$ and every term in that
identity is a power series.

For (ii), apply the same meromorphic identity to
$(X,Y,-X-Y)$ and write
\[
 \phi_1(T)=T^{-1}+F(T),\qquad
 \wp(T)=T^{-2}+W_0(T)
\]
for $T=X,Y,-X-Y$.  The purely polar terms cancel because, whenever
$X+Y+Z=0$,
\[
 (X^{-1}+Y^{-1}+Z^{-1})^2=X^{-2}+Y^{-2}+Z^{-2}.
\]
The remaining Laurent identity is \eqref{eq:all-integral-square}.
Oddness of $F$ gives
$P_0(0,Y)=P_0(X,0)=P_0(X,-X)=0$; hence
$XY(X+Y)\mid P_0$.  Finally,
\begin{align*}
 &F(X)F(Y)+F(Y)F(-X-Y)+F(-X-Y)F(X)\\
 &\qquad=\frac12\left(P_0^2-F(X)^2-F(Y)^2-F(-X-Y)^2\right),
\end{align*}
which turns \eqref{eq:all-integral-square} into
\eqref{eq:all-integral}.
\end{proof}

\subsection{Holomorphic projection of quadratic products}\label{sec:projection}

We apply algebraic holomorphic projection to the quadratic identities below.
Their images lie in the cuspidal coefficient quotient used in both the Manin
relation and the Fourier realization.

\begin{lemma}[Heat equation and chain rule]\label{lem:heat-chain}
Let $A(z,\tau)=(z-\bar z)/(\tau-\bar\tau)$.  Then
\begin{align}
 (\partial_z\log\theta(z,\tau))^2
 &=\phi_2(z,\tau)+\frac{2\pi i}{\tau-\bar\tau}
   +4\pi i\,\partial_\tau\log\theta(z,\tau),\label{eq:heat-square}\\
 \phi_1(z,\tau)^2
 &=\phi_2(z,\tau)+\frac{2\pi i}{\tau-\bar\tau}
   +4\pi iA(z,\tau)\phi_1(z,\tau)\notag\\
 &\qquad -(2\pi iA(z,\tau))^2
   +4\pi i\partial_\tau\log\theta(z,\tau).
 \label{eq:phi1-square}
\end{align}
For $z_0=a_1\tau+a_2$ and any smooth $f$,
\begin{equation}\label{eq:Taylor-chain}
 \Tay_{z_0}(\partial_\tau f)
 =\partial_\tau\Tay_{z_0}f-a_1\Tay_{z_0}(\partial_zf),
\end{equation}
where the derivative on the right differentiates the specialized
coefficients.
\end{lemma}

\begin{proof}
Differentiate \eqref{eq:phi1-theta} with respect to $z$ and use
$\partial_z\phi_1=-\phi_2$.  This gives
\[
 -\phi_2=\frac{2\pi i}{\tau-\bar\tau}
 +\partial_z^2\log\theta.
\]
The theta heat equation $\partial_z^2\theta=4\pi i\partial_\tau\theta$
rewrites the last term as
$4\pi i\partial_\tau\log\theta-(\partial_z\log\theta)^2$ and proves
\eqref{eq:heat-square}.  Substituting
$\partial_z\log\theta=\phi_1-2\pi iA$ gives
\eqref{eq:phi1-square}.  Finally, differentiate
$\partial_z^nf(z_0(\tau),\tau)$ by the chain rule, treating $\tau$ and
$\bar\tau$ as independent; since
$\partial_\tau z_0=a_1$ and $\partial_\tau\bar z_0=0$, summing the resulting
identity proves \eqref{eq:Taylor-chain}.
\end{proof}

\begin{proposition}[Quadratic decomposition]\label{prop:quadratic-decomp}
For $a\in\Q^2\setminus\Z^2$,
\begin{equation}\label{eq:quadratic-decomp}
 E_a(X)^2
 =E_{1,a}^2
 -2\Int_X\Tay_{\tau_a}\phi_3(\tau,X)
 +4\pi i\,\delta\Int_XE_a(X).
\end{equation}
\end{proposition}

\begin{proof}
Apply $\Tay_{\tau_a}$ to \eqref{eq:phi1-square}.  The elementary Taylor
identities
\begin{align*}
 \Tay_{\tau_a}\phi_2
 &=E_{2,a}-2\Int_X\Tay_{\tau_a}\phi_3,\\
 E_a&=E_{1,a}-\Int_X\Tay_{\tau_a}\phi_2,\\
 \Tay_{\tau_a}A&=a_1+\frac{X}{\tau-\bar\tau}
\end{align*}
follow from $\partial_z\phi_m=-m\phi_{m+1}$.  Integrating
\eqref{eq:phi1-theta} in the formal variable gives
\[
 \Tay_{\tau_a}\log\theta
 =\log\theta(\tau_a,\tau)+\Int_XE_a
 -2\pi ia_1X-\frac{\pi iX^2}{\tau-\bar\tau}.
\]
Use \eqref{eq:Taylor-chain} to compute the Taylor series of
$\partial_\tau\log\theta$.  After substitution in
\eqref{eq:phi1-square}, the constant term is $E_{1,a}^2$; the terms
containing $a_1$, $X/(\tau-\bar\tau)$, and
$X^2/(\tau-\bar\tau)^2$ cancel pairwise.  The remaining two series are
$-2\Int_X\Tay_{\tau_a}\phi_3$ and
$4\pi i\delta\Int_XE_a$, proving \eqref{eq:quadratic-decomp}.
\end{proof}

\begin{proposition}\label{prop:E1-square}
For $a\in\Q^2\setminus\Z^2$, the holomorphic modular form $E_{1,a}^2$
belongs to the Eisenstein subspace $\cE_2$ in the filtered colimit over
congruence subgroups.
\end{proposition}

\begin{proof}
Choose a level $L\ge5$ killing the class of $a$.  After an
$\SL(\Z/L\Z)$-change of coordinates, the transformation law in
Theorem~\ref{thm:torsion-specialization} reduces the assertion to a series
denoted $s_{b/L}$ by Borisov--Gunnells.  Their
\cite[Proposition~3.8]{BG1} states that $(s_{b/L})^2$ is Eisenstein.  The
normalization differs from $E_{1,a}$ by a nonzero scalar, so the same holds
for $E_{1,a}^2$.  Since this form is holomorphic, it lies in the holomorphic
Eisenstein subspace used here.
\end{proof}

\begin{theorem}[Projected square for a lattice coset]\label{thm:projected-square}
Let $V=a+r\Z^2$ be a single rational lattice coset, with
$a\in\Q^2$ and $r\in\Q^\times$.  Then
\[
 \cH\bigl(E_V(\tau,X)^2\bigr)=0\quad\text{in }\scrC_2.
\]
Equivalently, the holomorphic projection of $E_V^2$ has
Eisenstein-form coefficients.
\end{theorem}

\begin{proof}
Scaling reduces to $r=1$.  If the class of $a$ in $\Q^2/\Z^2$ is nonzero,
then $E_V=E_a$.  In \eqref{eq:quadratic-decomp}, the constant is Eisenstein
by Proposition~\ref{prop:E1-square}; every coefficient of the integrated
$\phi_3$ term is an Eisenstein series by
Theorem~\ref{thm:torsion-specialization}; and algebraic holomorphic projection
kills the raised term by Theorem~\ref{thm:H-package}.  If $a\in\Z^2$, then
$E_V=F$.  Formula \eqref{eq:F-square} has the same three types of terms:
$-\partial_XF-2F/X-3E_{2,0}$ has Eisenstein coefficients and the remaining
term is raised.  Its image in $\scrC_2$ therefore vanishes.
\end{proof}

\begin{remark}
The assertion is for a single coset; the square of a sum of coset series
also contains cross terms.
\end{remark}

\section{The modular symbol}\label{sec:symbol}

\subsection{The distribution on \texorpdfstring{$\Mat(\Q)$}{M2(Q)}}\label{sec:distribution}

For $\varphi\in\cS(\Q^2)$, let $E(\varphi;X)$ denote the linear extension
of $[V]\mapsto E_V(X)$.

\begin{construction}\label{constr:mu}
Define the holomorphic distribution and its cuspidal quotient by
\[
 \widetilde\mu\in\Dist\bigl(\Mat(\Q),\scrM_2\bigr),
 \qquad
 \mu\in\Dist\bigl(\Mat(\Q),\scrC_2\bigr)
\]
on decomposable test functions by
\begin{equation}\label{eq:mu-test}
 \begin{aligned}
 \widetilde\mu(\varphi_1\otimes\varphi_2)
 &=\widetilde{\cH}\bigl(E(\varphi_1;X)E(\varphi_2;Y)\bigr),\\
 \mu(\varphi_1\otimes\varphi_2)
 &=\cH\bigl(E(\varphi_1;X)E(\varphi_2;Y)\bigr).
 \end{aligned}
\end{equation}
Thus $\mu$ is the coefficientwise quotient of $\widetilde\mu$.  Via
\eqref{eq:tensor-test}, both distributions are defined on all test functions.
\end{construction}

\begin{proposition}[Well-definedness and distribution relation]\label{prop:mu-well-defined}
The two assignments in \eqref{eq:mu-test} are bilinear and therefore define
unique distributions.  For every finite disjoint decomposition
$U=\coprod_iU_i$ with $[U_i]\in\cS(\Mat(\Q))$,
\[
 \mu(U)=\sum_i\mu(U_i).
\]
If $U=V_1\times V_2$ under the column identification $\Mat(\Q)\simeq\Q^2\times\Q^2$, then
\[
 \mu(U)=\cH\bigl(E_{V_1}(X)E_{V_2}(Y)\bigr).
\]
\end{proposition}

\begin{proof}
The map $\varphi\mapsto E(\varphi;X)$ is linear by
Proposition~\ref{prop:compact-cont}; multiplication of power series is
bilinear, and $\widetilde{\cH}$ and $\cH$ are linear.  Thus the two assignments in
\eqref{eq:mu-test} are $\Z$-bilinear on
$\cS(\Q^2)\times\cS(\Q^2)$.  The universal property of the tensor product and
\eqref{eq:tensor-test} give unique homomorphisms on
$\cS(\Mat(\Q))$.  Applying $\mu$ to
$[U]=\sum_i[U_i]$ proves finite additivity.  The last formula is simply
\eqref{eq:mu-test} applied to $[V_1]\otimes[V_2]$.
\end{proof}

\subsubsection{Torus invariance}\label{sec:torus}

Let $T(\Q)\subset\G(\Q)$ be the diagonal torus.

\begin{proposition}[Torus invariance]\label{prop:torus}
For every $t\in T(\Q)$,
\[
 t\cdot\mu=\mu.
\]
\end{proposition}

\begin{proof}
Write $t=\operatorname{diag}(a,d)$.  It is enough to evaluate on a
rectangle $V_1\times V_2$, since their characteristic functions generate the
$\Z$-module $\cS(\Mat(\Q))$.  Then
\begin{align*}
 (t\cdot\mu)(V_1\times V_2)(X,Y)
 &=ad\,\mu(aV_1\times dV_2)(aX,dY)\\
 &=ad\,\cH\bigl(E_{aV_1}(aX)E_{dV_2}(dY)\bigr)\\
 &=\cH\bigl(E_{V_1}(X)E_{V_2}(Y)\bigr),
\end{align*}
by Proposition~\ref{prop:scaling}, which includes the regularized zero coset.
\end{proof}

\subsection{Manin relations and the Eisenstein--Kronecker modular symbol}\label{sec:Manin}

The two Manin relations determine $\Phi$ from its value $\Phi(D_\infty)=\mu$.

\begin{theorem}[Manin presentation]\label{thm:Manin-presentation}
As a left $\Z[\G(\Q)]$-module, $\Delta_0$ is generated by $D_\infty$.  Its annihilator is the left ideal generated by
\[
 \{t-1:t\in T(\Q)\},\qquad 1+S,\qquad 1+R+R^2.
\]
Consequently, for a left $\G(\Q)$-module $W$, evaluation at $D_\infty$ identifies $\Symb_{\G(\Q)}(W)$ with the elements $w\in W^{T(\Q)}$ satisfying
\[
 (1+S)w=0,
 \qquad
 (1+R+R^2)w=0.
\]
\end{theorem}

This is the rank-two Manin presentation; compare
\cite[Proposition~3.3.1 and Theorem~1.1.4]{BPPS} and \cite{Stevens}.  We include
the short proof in order to fix the left-action, orientation, and torus
conventions used below.

\begin{proof}
For distinct $r,s\in\mathbf P^1(\Q)$ write $[r,s]=[r]-[s]$.
Every such ordered pair is $\gamma[\infty,0]$ for some
$\gamma\in\G(\Q)$; two choices differ on the right by an element of the
diagonal torus.  The relations $[r,s]+[s,r]=0$ and
$[r,s]+[s,t]+[t,r]=0$ generate all relations among degree-zero divisors:
fixing a base cusp $s_0$, one has $[r,s]=[r,s_0]-[s,s_0]$.  The first
relation is the translate of $(1+S)D_\infty=0$, and the relation for the
ordered triple $(\infty,0,1)$ is
$(1+R+R^2)D_\infty=0$.  Since $\operatorname{PGL}_2(\Q)$ is sharply
three-transitive on ordered triples of distinct rational points, every
triangle relation is its translate.  This proves the asserted presentation.
Applying $\operatorname{Hom}_{\G(\Q)}(-,W)$ gives the last statement.
\end{proof}

\begin{proposition}[First Manin relation]\label{prop:Manin1}
The distribution $\mu$ satisfies
\[
 (1+S)\mu=0.
\]
\end{proposition}

\begin{proof}
It is enough to evaluate on a lattice coset in $\Mat(\Q)$
$U=A+r\Mat(\Z)$.  Their characteristic functions generate
$\cS(\Mat(\Q))$ as a $\Z$-module.  Right
multiplication by $S$ exchanges its two column cosets and changes the sign of
one of them.  Proposition~\ref{prop:scaling}, parity, and commutativity give
\[
 (S\mu)(U)(X,Y)
 =\cH\bigl(E_{U_2}(Y)E_{-U_1}(-X)\bigr)
 =-\mu(U)(X,Y).
\]
\end{proof}

\begin{proposition}[Second Manin relation]\label{prop:Manin2}
The distribution $\mu$ satisfies
\[
 (1+R+R^2)\mu=0.
\]
\end{proposition}

\begin{proof}
It suffices to evaluate on $U=A_0+r\Mat(\Z)$.  Scalar-torus invariance and
Proposition~\ref{prop:scaling} reduce the calculation to $r=1$.  Let
$a,b\in\Q^2/\Z^2$ be the two column classes of $A_0$ and put
$c=-a-b$.  Right multiplication by $1,R,R^2$ gives the three products
\[
 AB,\qquad BC,\qquad CA,\qquad
 A=E_a(X),\ B=E_b(Y),\ C=E_c(-X-Y),
\]
where $E_0$ means the regularized series $F$.  Thus the desired value is
$\cH(AB+BC+CA)$, and
\begin{equation}\label{eq:ABC}
 AB+BC+CA=\frac12(A+B+C)^2-\frac12(A^2+B^2+C^2).
\end{equation}
Since $a+b+c=0$, either none, one, or all three parameters are zero.

\smallskip\noindent
\emph{No parameter is zero.}
The fundamental identity, Theorem~\ref{thm:fundamental-identity}, gives
\[
 (A+B+C)^2=W_a(X)+W_b(Y)+W_c(-X-Y),
 \qquad W_u(T)=\Tay_{\tau_u}\wp(T).
\]
Every coefficient on the right is Eisenstein by
Theorem~\ref{thm:torsion-specialization}.  Each square in
\eqref{eq:ABC} has Eisenstein holomorphic projection by
Theorem~\ref{thm:projected-square}; hence the cuspidal image is zero.

\smallskip\noindent
\emph{Exactly one parameter is zero.}
After a cyclic permutation assume $b=0$, so $c=-a\ne0$ and $B=F(Y)$.
Put $P=A+B+C=P_a(X,Y)$.  The regularized fundamental identity
\eqref{eq:one-integral} applies.  Proposition~\ref{prop:regularized-fundamental}
also shows that $P/Y$ is a power series.  Its
coefficients are finite differences of the Eisenstein coefficients of
$E_{-a}$ together with coefficients of $F/Y$.  The only
almost-holomorphic coefficient of $F/Y$ is its constant coefficient
$-E_{2,0}$.  In $P/Y$ it occurs as the difference
$E_{2,-a}-E_{2,0}$, which is holomorphic and Eisenstein by
Theorem~\ref{thm:torsion-specialization}.  Hence every coefficient of
$P/Y$ is holomorphic Eisenstein.
The first three terms in \eqref{eq:one-integral} are also Eisenstein, and the
three individual squares are handled by
Theorem~\ref{thm:projected-square}.  Equation \eqref{eq:ABC} again has zero
cuspidal image.

\smallskip\noindent
\emph{All three parameters are zero.}
Put
\[
 P=P_0(X,Y)=F(X)+F(Y)+F(-X-Y),
 \qquad Q=\frac1X+\frac1Y-\frac1{X+Y}.
\]
The fully regularized identity is \eqref{eq:all-integral}.  By
Proposition~\ref{prop:regularized-fundamental}, $PQ$ is a power series.  Its coefficients are
linear combinations of the level-one Eisenstein series occurring in $F$.
The possible nonholomorphic term cancels already in $P$: its contribution is
a scalar multiple of $E_{2,0}(X+Y-X-Y)=0$.
The remaining terms in \eqref{eq:all-integral} are Eisenstein after
holomorphic projection by Proposition~\ref{prop:regularized}.  Thus the
cuspidal image is zero in the last case as well.

We have proved the identity on a generating set of test functions, hence
as an identity of distributions.
\end{proof}

\begin{proof}[Proof of Theorem~\ref{thm:intro-main}]
By Proposition~\ref{prop:torus}, $\mu$ is $T(\Q)$-invariant.  Propositions~\ref{prop:Manin1} and \ref{prop:Manin2} give the two Manin relations.  Theorem~\ref{thm:Manin-presentation} therefore gives a unique $\G(\Q)$-modular symbol $\Phi$ such that $\Phi(D_\infty)=\mu$.
\end{proof}

\subsection{Specialization to Borisov--Gunnells symbols}\label{sec:BG}

We now give the finite-level comparison with Borisov--Gunnells.

Let
\[
 U_1(N)=
 \begin{pmatrix}0&0\\0&1/N\end{pmatrix}
 +\Mat(\Z).
\]
\begin{definition}
For $k\geq2$, define
\[
 \Phi_{k,N}(D)
 =\Phi(D)(U_1(N))^{[k-2]}.
\]
\end{definition}

\begin{proposition}[Level and weight]\label{prop:level-weight}
The map $\Phi_{k,N}$ takes values in
\[
 \cC_k(\Gamma_1(N))\otimes\operatorname{Sym}^{k-2}(\C^2)
\]
and is $\Gamma_1(N)$-equivariant.  Under the cuspidal splitting, it is cusp-form-valued.
\end{proposition}

\begin{proof}
The right action of $\Gamma_1(N)$ on $\Mat(\Q)/\Mat(\Z)$ fixes
$U_1(N)$.  Indeed, for
$\gamma=\begin{psmallmatrix}a&b\\ c&d\end{psmallmatrix}\in\Gamma_1(N)$ one has
$c\equiv0\pmod N$ and $d\equiv1\pmod N$, and hence
\[
 \begin{pmatrix}0&0\\0&1/N\end{pmatrix}\gamma
 =\begin{pmatrix}0&0\\c/N&d/N\end{pmatrix}
 \equiv
 \begin{pmatrix}0&0\\0&1/N\end{pmatrix}
 \pmod{\Mat(\Z)}.
\]
Thus evaluation of the distribution $\Phi(D)$ on $U_1(N)$ is compatible with
restriction from $\SL(\Z)$ to $\Gamma_1(N)$.  The equivariance of $\Phi$ then
gives the usual modular-symbol transformation law for $\Phi_{k,N}$.  Taking
the homogeneous degree $k-2$ part gives the coefficient representation
$\operatorname{Sym}^{k-2}(\C^2)$, while the holomorphic projection places the
modular-form coefficient in weight $k$ and level $\Gamma_1(N)$, modulo the
Eisenstein part.  After choosing the cuspidal splitting of the coefficient
module, this is represented by cusp forms.
\end{proof}

\begin{theorem}[Explicit Borisov--Gunnells comparison]\label{thm:BG-explicit}
For $\gamma=\begin{psmallmatrix}a&b\\c&d\end{psmallmatrix}\in\SL(\Z)$, the Manin-symbol value of $\Phi_{k,N}$ at the class of $\gamma$ is the homogeneous degree-$k-2$ part of
\[
 \cH\left(
 E_{(0,c/N)}\bigl(aX+cY\bigr)
 E_{(0,d/N)}\bigl(bX+dY\bigr)
 \right),
\]
with regularized interpretation if $c$ or $d$ is zero modulo $N$.  More
explicitly, its degree-$k-2$ representative is
\begin{equation}\label{eq:BG-coeff}
 (-1)^{k-2}\prcusp H_k^{\mathrm{alg}}
 \sum_{r=0}^{k-2}E_{r+1,(0,c/N)}E_{k-r-1,(0,d/N)}
 (aX+cY)^r(bX+dY)^{k-2-r}.
\end{equation}
Applying finite Fourier transform to the complete family of coefficient
polynomials indexed by $(c,d)\in(\Z/N\Z)^2$, and then the pairing
\eqref{eq:BG-pairing}, gives \cite[Equation~(8), \S6]{BG2} after
holomorphic projection and passage modulo Eisenstein series.
The restriction to primitive pairs is made after the Fourier transform.
\end{theorem}

\begin{proof}
The right action of
$\gamma=\begin{psmallmatrix}a&b\\c&d\end{psmallmatrix}$ sends the coset
$U_1(N)$ to a lattice coset whose column classes are $(0,c/N)$ and $(0,d/N)$.
The coefficient action sends $(X,Y)$ to
$(aX+cY,bX+dY)$.  Therefore equivariance of $\Phi$ gives the displayed
generating series.  Taking homogeneous degree $k-2$ and using
\eqref{eq:Ea} gives \eqref{eq:BG-coeff}; this proves the formula in every
weight, including the signs and the two endpoint terms.  If a parameter is
zero, the omitted-term regularization applies; in particular, $E_{1,0}=0$.

Set
\[
 s_t^{(m)}=E_{m,(0,t/N)},\qquad
 C_m=\frac{(m-1)!}{(-2\pi i)^m},\qquad t\in\Z/N\Z.
\]
Define the completed residue-class series by
\begin{equation}\label{eq:BG-fourier}
 \widehat s_u^{(m)}
 =\frac{C_m}{N}\sum_{t\bmod N}e^{-2\pi iut/N}s_t^{(m)}.
\end{equation}
For $0\le u<N$, their Fourier expansions are
\[
 \widehat s_u^{(m)}
 =b_{m,u}+\sum_{n\ge1}q^n\sum_{d\mid n}d^{m-1}
   \bigl(\one_{d\equiv u\, (N)}+(-1)^m\one_{d\equiv-u\, (N)}\bigr)
   +\frac{\one_{m=2}\one_{u=0}}{4\pi y},
\]
where
\[
 b_{m,u}=-\frac{N^{m-1}}m B_m(u/N),
 \qquad b_{1,0}=0.
\]
Here $B_m(T)$ is the Bernoulli polynomial.  The positive coefficients
follow from character orthogonality.  For $m\ge2$, the constant term
comes from
$C_mN^{m-1}\sum_{\ell\ne0}e^{-2\pi iu\ell/N}\ell^{-m}$;
the Fourier series of $B_m$ gives the displayed value.  In weight one
the finite transform of the cotangent gives $1/2-u/N$ for $u\ne0$,
and parity gives zero at $u=0$.  The common nonholomorphic term of the
$C_2s_t^{(2)}$ survives the transform only at $u=0$.
Thus $\widehat s_u^{(m)}$ is the generator $\widetilde s_{u/N}^{(m)}$
of \cite[\S2]{BG2}, completed when $(m,u)=(2,0)$; the zero weight-one
generator is interpreted as zero.

To specify the polynomial normalization, put $n=k-2$ and
\[
 \begin{aligned}
 &P_{c,d}(X,Y)=(-1)^n\sum_{a+b=n}s_c^{(a+1)}s_d^{(b+1)}X^aY^b,\\
 &\widehat P_{u,v}=\frac1{N^2}\sum_{c,d\bmod N}
   e^{-2\pi i(uc+vd)/N}P_{c,d}.
 \end{aligned}
\]
These are the coefficient polynomials underlying \eqref{eq:BG-coeff},
before substitution by $\gamma$ and before projection.  For a homogeneous
polynomial $p(x,y)$ of degree $n$, define
\begin{equation}\label{eq:BG-pairing}
 \langle p,P\rangle
 =\left.p(\partial_Y,-\partial_X)P(X,Y)\right|_{X=Y=0}.
\end{equation}
The pairing is invariant under simultaneous $\SL$-substitutions in
the two polynomials, by the chain rule.
Equation~\eqref{eq:BG-fourier} gives
\[
 \widehat P_{u,v}
 =(-1)^n(-2\pi i)^k
   \sum_{a+b=n}\widehat s_u^{(a+1)}\widehat s_v^{(b+1)}
                    \frac{X^aY^b}{a!b!}.
\]
Consequently, for $r+s=n$,
\begin{equation}\label{eq:BG-exact-comparison}
 \frac{\langle x^ry^s,\widehat P_{u,v}\rangle}
      {(-1)^n(-2\pi i)^k}
 =(-1)^s\widehat s_u^{(s+1)}\widehat s_v^{(r+1)}.
\end{equation}
This is the completed form of \cite[Equation~(8)]{BG2}.

We verify compatibility with the Manin relations before restricting the
indices.  Write $t=(c,d)$ and $z=(X,Y)$.  If a matrix $g\in\SL(\Z)$
acts on a polynomial family by $P_t(z)\mapsto P_{tg}(zg)$, change of
summation variable gives
\[
 \widehat{(g\cdot P)}_u(z)=\widehat P_{u g^{-t}}(zg),\qquad
 \langle p,\widehat{(g\cdot P)}_u\rangle
 =\langle p\circ g^{-1},\widehat P_{u g^{-t}}\rangle.
\]
Here $g^{-t}$ is the inverse transpose and $(p\circ g)(z)=p(zg)$.
For the projected family $\cH(P)$, the $S$-relation holds, as does the
cyclic relation for $R^{-t}=SRS^{-1}$: this follows from the $R$-relation
and $S^{-1}\mu=-\mu$.  Taking $g=S$ and $g=R^{-t}$ in the displayed
identities therefore gives precisely the relations of
\cite[\S3.1]{BG2}, in which the index transforms by $S,R$ and the
polynomial by $S^t,R^t$.  Both matrices preserve primitive pairs.
Thus \eqref{eq:BG-exact-comparison}, after projection and restriction,
is an equality of maps on Manin symbols, not only of coefficient formulas.

For $k>2$ this identifies the projected map of
\cite[Theorem~6.2 and Corollary~6.3]{BG2}; its derivative terms complete
to raised Eisenstein series and are killed by the projection.  For $k=2$,
$E_{1,(0,t/N)}/(2\pi i)$ is the theta-derivative series of
\cite[Definition~3.11]{BG1}.  Hence the untransformed specialization,
divided by $(2\pi i)^2$, is their weight-two map.  The Fourier-transformed
formula is its Fricke transform, by \cite[Proposition~3.19]{BG1}.
The argument applies also at levels $2,3,4$: the Fourier identity and the
two Manin relations hold on every pair modulo $N$, including zero.
\end{proof}

\begin{proof}[Proof of Theorem~\ref{thm:intro-BG}]
The coset $U_1(N)$ is stable under right multiplication by $\Gamma_1(N)$, so
equivariance of $\Phi$ and homogeneity give a $\Gamma_1(N)$-modular symbol.
Theorem~\ref{thm:BG-explicit} recovers the Borisov--Gunnells map from
the full finite family by \eqref{eq:BG-exact-comparison}.
The quotient $M_k/\cE_k$ is canonically the cuspidal quotient; applying the
chosen cusp projection gives the stated cusp-form-valued representative.
\end{proof}

\section{Finite Fourier transform and uniqueness}\label{sec:Colmez}

The arithmetic comparison involves products of two different Eisenstein
families.  We first relate these families by finite Fourier transform,
keeping track of the normalization and of rational changes of scale.
We then prove that the torus and Manin relations determine a distribution
from the partial Fourier values needed in Section~\ref{sec:BK}.

We use the two Eisenstein families of \cite[\S\S1.3.1--1.3.2]{Colmez}, with their
full-lattice normalization.  For $x\in(\Q/\Z)^2$, put
\begin{equation}\label{eq:Colmez-Pasol-normalization}
 \begin{aligned}
 &C_k=\frac{(k-1)!}{(-2\pi i)^k},\qquad
 E_x^{(k)}=C_kE_{k,x},\\
 &F_x^{(k)}=C_k\left.
 \sum_{(u,v)\in\Z^2\setminus\{(0,0)\}}
 \frac{e^{2\pi i(ux_2-vx_1)}}{(u\tau+v)^k|u\tau+v|^s}
 \right|_{s=0},
 \end{aligned}
\end{equation}
where the last value is defined by analytic continuation.  The first
family translates the denominator; the second inserts an additive
character in the numerator.  In particular,
\begin{equation}\label{eq:completed-zero-normalization}
 E_0^{(2)}=F_0^{(2)}=-\frac1{12}E_2^*.
\end{equation}
Indeed, the expansion
$\partial_z\log\theta(z,\tau)=z^{-1}-\frac{\pi^2}{3}E_2(\tau)z+O(z^3)$
in \eqref{eq:phi1-theta} gives
$E_{2,0}=\frac{\pi^2}{3}E_2^*$, and $C_2=-1/(4\pi^2)$.
Here $E_2=E_2^*+3/(\pi\operatorname{Im}\tau)$.
The value \eqref{eq:completed-zero-normalization} follows from the defining
full-lattice sums; it differs by a factor of two from the displayed value
in \cite[Proposition~1.1(i)]{Colmez} and
\cite[Proposition~2.8(1)]{Wang}.

For $x=(x_1,x_2)$ and $y=(y_1,y_2)$ in $(\Z/N\Z)^2$, put
\[
 \langle x,y\rangle=x_1y_2-x_2y_1\in\Z/N\Z.
\]
In a finite-level family, an index $x\in(\Z/N\Z)^2$ denotes the torsion
parameter $x/N$.
Let $\scrS_N$ denote the unnormalized symplectic Fourier transform
\begin{equation}\label{eq:finite-FT}
 (\scrS_Nf)(x)
 =\sum_{y\in(\Z/N\Z)^2}
 \exp\left(\frac{2\pi i}{N}\langle x,y\rangle\right)f(y),
\end{equation}
and put
\begin{equation}\label{eq:normalized-finite-FT}
 \cF_{N,k}=(-1)^kN^{-k}\scrS_N.
\end{equation}

Brunault and Zudilin use the finite Fourier transform between the two
torsion families to derive a weight-three product relation
\cite[\S6, proof of Theorem~56]{BrunaultZudilin}.  For the comparison
with Colmez's distribution we record the transform in the normalization
of \eqref{eq:Colmez-Pasol-normalization}, including the zero parameter.

\begin{theorem}[Finite Fourier transform]\label{thm:finite-Fourier}
With the normalization of \eqref{eq:finite-FT}, the torsion families satisfy
\[
 \cF_{N,k}\bigl(E^{(k)}\bigr)=F^{(k)},
 \qquad
 E^{(k)}=(-1)^kN^{k-2}\scrS_N\bigl(F^{(k)}\bigr).
\]
Equivalently,
\begin{equation}\label{eq:unnormalized-Fourier-EF}
 \scrS_NE^{(k)}=(-1)^kN^kF^{(k)},
 \qquad
 \scrS_NF^{(k)}=(-1)^kN^{2-k}E^{(k)}.
\end{equation}
The identities hold also at the zero torsion parameter when the omitted-term
regularization is used, and they are compatible with the
$\SL(\Z/N\Z)$-action.
\end{theorem}

\begin{proof}
For $\operatorname{Re}(s)$ in the half-plane of absolute convergence,
decompose the lattice into its $N^2$ residue classes and apply character
orthogonality.  With the positive kernel in \eqref{eq:finite-FT}, the residue
calculation gives the factor $N^{k+s}$ and the character indexed by $-x$.
Continuation to $s=0$ therefore gives
\[
 \scrS_NE^{(k)}(x)=N^kF^{(k)}(-x)
 =(-1)^kN^kF^{(k)}(x).
\]
Both continued families are regular at this special value.  At $x=0$ the
same argument means that the singular lattice term is omitted on both sides.
Since character orthogonality also gives $\scrS_N^2=N^2\operatorname{id}$,
the second formula in \eqref{eq:unnormalized-Fourier-EF} follows.  Invariance of the
alternating pairing proves $\SL(\Z/N\Z)$-equivariance.
\end{proof}

\begin{theorem}[Colmez's Eisenstein distributions]\label{thm:Colmez-Eis}
For $m\geq1$ there are algebraic distributions
$z_{\Eis}(m)$ and $z'_{\Eis}(m)$ on $\cS(\Q^2)$ with values in
weight-$m$ Eisenstein series such that, for $r\in\Q^\times$,
\begin{align*}
 \int_{(a+r\Z)\times(b+r\Z)}z_{\Eis}(m)
 &=r^{-m}E^{(m)}_{r^{-1}a,r^{-1}b},\\
 \int_{(a+r\Z)\times(b+r\Z)}z'_{\Eis}(m)
 &=r^{m-2}F^{(m)}_{r^{-1}a,r^{-1}b}.
\end{align*}
At weight two the values are the completed, omitted-term series of
\eqref{eq:Colmez-Pasol-normalization}.
Their determinant covariance is given by \cite[Theorem~1.5]{Colmez}.
\end{theorem}

\begin{definition}[Colmez]\label{def:Colmez-product}
For $k\geq2$ and $1\leq j\leq k-1$, Colmez's product distribution is
\[
 z_{\Eis}(k,j)
 =\frac{(-1)^j}{(j-1)!}\,
 z'_{\Eis}(k-j)\otimes z_{\Eis}(j),
\]
where $z'_{\Eis}$ has values in the $F$-family and $z_{\Eis}$ in the
$E$-family \cite[\S1.3.5]{Colmez}.
\end{definition}

Put $m=k-j$, and let $\mu_k^{(j)}$ denote the coefficient of
$X^{m-1}Y^{j-1}$ in the weight-$k$ part of $\mu$, and set
\begin{equation}\label{eq:comparison-transform}
 \mathscr F_{N,k,j}
 =\frac{(-2\pi i)^k}{(m-1)!}N^{m-2}\scrS_N^{(1)} .
\end{equation}
Here $\scrS_N^{(1)}$ acts on the first row.

\begin{theorem}[Finite-level coefficient identity]\label{thm:Colmez-coeff}
At level $N$ one has
\begin{equation}\label{eq:Colmez-comparison}
 \mu_k^{(j)}
 =\pi_{\cC_k}\,H_k^{\mathrm{alg}}
 \bigl(\mathscr F_{N,k,j}z_{\Eis}(k,j)\bigr).
\end{equation}
The identity includes the zero torsion parameter by omitted-term
regularization.  In weight two it is used after $H_k^{\mathrm{alg}}$ and
$\pi_{\cC_k}$.
\end{theorem}

\begin{proof}
At a torsion pair $(a,b)$ the finite-level function underlying Colmez's
product distribution is
\[
 \frac{(-1)^j}{(j-1)!}F_a^{(m)}E_b^{(j)}.
\]
Applying \eqref{eq:unnormalized-Fourier-EF} in the first row gives
\begin{align*}
 \scrS_N^{(1)}\left(
  \frac{(-1)^j}{(j-1)!}F^{(m)}E^{(j)}\right)(a,b)
 &=\frac{(-1)^{j+m}N^{2-m}}{(j-1)!}
   E_a^{(m)}E_b^{(j)}\\
 &=(-1)^kN^{2-m}\frac{(m-1)!}{(-2\pi i)^k}
   E_{m,a}E_{j,b}.
\end{align*}
Equation~\eqref{eq:Ea} gives
\[
 [X^{m-1}Y^{j-1}]E_a(X)E_b(Y)
 =(-1)^{m-1+j-1}E_{m,a}E_{j,b}
 =(-1)^kE_{m,a}E_{j,b}.
\]
Multiplication by the scalar in \eqref{eq:comparison-transform}, followed by
$H_k^{\mathrm{alg}}$ and the cuspidal quotient, proves
\eqref{eq:Colmez-comparison}.

The calculation remains valid when a parameter is zero.  Odd weights
vanish there by parity.  In even weight the omitted-term continuation is
the regularized level-one series, with the completed weight-two value
\eqref{eq:completed-zero-normalization}.  Thus the product identity holds
in the nearly holomorphic coefficient space before projection.  Applying
the linear maps $H_k^{\mathrm{alg}}$ and $\pi_{\cC_k}$ preserves this
identity; no multiplicativity of holomorphic projection is used.
\end{proof}

For later use, fix $r\in\Q_{>0}$ and $x,y\in N^{-1}\Z^2/\Z^2$, and write
\begin{equation}\label{eq:adelic-row-coset}
 B_{r,N}(x,y)
 =\bigl(rx+r\Z^2\bigr)\times
   \bigl(ry+r\Z^2\bigr)
 \subset\Mat(\Q),
\end{equation}
where the two factors are the two rows.  To compare the values on these
cosets, the operator \eqref{eq:comparison-transform} must be multiplied
by the scale factor $r^{2-2m}$.  Thus put
\begin{equation}\label{eq:scaled-comparison-transform}
 \mathscr F^{(r)}_{N,k,j}
 =r^{2-2m}\mathscr F_{N,k,j}
 =\frac{(-2\pi i)^k}{(m-1)!}
   r^{2-2m}N^{m-2}\scrS_N^{(1)}.
\end{equation}

\begin{corollary}[Rational lattice cosets]
\label{cor:scaled-coset-comparison}
On \eqref{eq:adelic-row-coset}, before holomorphic projection one has
\begin{equation}\label{eq:scaled-coset-value}
 \mathscr F^{(r)}_{N,k,j}z_{\Eis}(k,j)
 =(-1)^kr^{-k}E_{m,x}E_{j,y}.
\end{equation}
After applying $H_k^{\mathrm{alg}}$ and $\pi_{\cC_k}$ one obtains the value of $\mu_k^{(j)}$ on the corresponding transposed lattice coset in $\Mat(\Q)$.
\end{corollary}

\begin{proof}
Colmez's formulas for compact open subsets give
\[
 z_{\Eis}(k,j)\bigl(B_{r,N}(x,y)\bigr)
 =\frac{(-1)^j}{(j-1)!}
   r^{m-j-2}F_x^{(m)}E_y^{(j)}.
\]
The Fourier calculation before the scale factor in the proof of
Theorem~\ref{thm:Colmez-coeff} contributes
$(-1)^mN^{2-m}$ and changes the first $F$-factor into an $E$-factor.
The factor $r^{2-2m}$ in
\eqref{eq:scaled-comparison-transform} changes the scale exponent by
\[
 (m-j-2)+(2-2m)=-m-j=-k.
\]
The factorial and period calculation is unchanged, proving
\eqref{eq:scaled-coset-value}.  In the column convention, the coefficient of
the series attached to a coset of scale $r$ is also
$(-1)^kr^{-m}r^{-j}E_{m,x}E_{j,y}$, by Proposition~\ref{prop:scaling}.
Transposition and coefficientwise projection give the last assertion.
\end{proof}

\subsection{Uniqueness from finite Fourier values}

Right multiplication by the lower unipotent matrix
$L=\left(\begin{smallmatrix}1&0\\1&1\end{smallmatrix}\right)$
adds the second column of a box to the first.  The finite Fourier
transform diagonalizes this translation.  We use its nontrivial
eigenspaces to determine a distribution satisfying the Manin relations.

Put $T=\cS(\Mat(\Q))\otimes_\Z\C$.  Retain the tests
$e^{r,N}_{a,b}$ and $f^{r,N}_{x,b}$ of the introduction, and put
\[
 T_{r,N}=\operatorname{span}_{\C}\{e^{r,N}_{a,b}\},\qquad
 A_{r,N}=\operatorname{span}_{\C}
 \{f^{r,N}_{x,b}:\langle x,b\rangle\not\equiv0\pmod N\}.
\]
Every rational box occurs in some $T_{r,N}$, and these spaces span $T$.
The condition defining $A_{r,N}$ is equivalent to
$e^{2\pi i\langle x,b\rangle/N}\ne1$, also for composite $N$.

\begin{theorem}\label{thm:Fourier-Manin-uniqueness}
Let $V$ be a complex vector space, $n\ge0$, and
\[
 \beta:T\longrightarrow V\otimes_\C\C[X,Y]_n,\qquad
 \beta(U)=\sum_{i=0}^na_i(U)X^{n-i}Y^i.
\]
Give such maps the action \eqref{eq:old-action}, acting trivially on $V$.
Suppose that $\beta$ is invariant under the rational diagonal torus and
satisfies $(1+S)\beta=(1+R+R^2)\beta=0$.  Suppose that $a_0$ vanishes
on every $A_{r,N}$ and, if $n\ge2$, that $a_n$ also vanishes on every
$A_{r,N}$.  Then $\beta=0$.
\end{theorem}

\begin{proof}
Let $\tau[U]=[UL]$.  Then
\[
 \tau e_{a,b}=e_{a+b,b},\qquad
 \tau f_{x,b}=e^{-2\pi i\langle x,b\rangle/N}f_{x,b}.
\]
Character orthogonality gives a basis $\{f_{x,b}\}$ of $T_{r,N}$; hence
\begin{equation}\label{eq:shear-image}
 (\tau-1)T_{r,N}=A_{r,N}.
\end{equation}
Thus each coefficient assumed to vanish on every $A_{r,N}$ satisfies
$a_i(UL)=a_i(U)$.  Torus invariance gives, for every $i$,
\begin{equation}\label{eq:uniqueness-torus}
 a_i(U\operatorname{diag}(u,v))
 =u^{-(n-i+1)}v^{-(i+1)}a_i(U).
\end{equation}
Conjugating $L$ by diagonal matrices proves lower-unipotent invariance
of $a_0$ and, when $n\ge2$, of $a_n$ over $\Q$.

Suppose $n\ge2$.  The $S$-relation gives
\[
 a_0(U)+(-1)^na_n(US)=0,\qquad a_n(U)+a_0(US)=0.
\]
Thus $a_0,a_n$ are also upper-unipotent invariant, hence
$\SL(\Q)$-invariant.  Indeed, elementary elimination writes an element
of $\SL(\Q)$ as a product of upper and lower unipotents.  Taking
$(u,v)=(2,1/2)$ in \eqref{eq:uniqueness-torus} gives
$a_0=2^{-n}a_0$ and $a_n=2^na_n$, so both vanish.  The coefficient
of $X^n$ in the cyclic relation is now
\[
 \sum_{i=1}^{n-1}(-1)^{n-i}a_i(UR^2)=0.
\]
Replace $UR^2$ by $U\operatorname{diag}(t,1)$ and use
\eqref{eq:uniqueness-torus}.  Distinct powers of $t^{-1}$ have zero
coefficients.  Thus every $a_i$ vanishes.

Suppose $n=1$.  Then $a_1(U)=-a_0(US)$.  Since $RS=-L$,
\[
 a_1(UR)=-a_0(-UL)=a_0(U).
\]
The coefficient of $X$ in the cyclic relation is
\[
 a_0(U)-a_1(UR)-a_0(UR^2)+a_1(UR^2)=0.
\]
Consequently $a_1=a_0$.  Under $U\mapsto U\operatorname{diag}(2,1)$,
these two functions scale respectively by $1/2$ and $1/4$.
Both vanish.

If $n=0$, the $S$-relation again gives upper-unipotent invariance from
lower-unipotent invariance.  Therefore $a_0(US)=a_0(U)$, whereas
$a_0(US)=-a_0(U)$.  Thus $a_0=0$.
\end{proof}

\begin{proof}[Proof of Corollary~\ref{cor:endpoint-reconstruction}]
Fix $k$ and take the difference of the values at $D_\infty$ of two
weight-$k$ symbols with the same endpoint data.  This difference
satisfies the torus and Manin relations.  Put $n=k-2$ in
Theorem~\ref{thm:Fourier-Manin-uniqueness}.  The index $j=1$ gives
$a_0$, and for $k\ge4$ the index $j=k-1$ gives $a_n$.  Agreement on
the stated tests means that these coefficients vanish on every
$A_{r,N}$.  The theorem gives equality at $D_\infty$, including in
degrees $n=0,1$ where $a_0$ alone is tested.  Equality of the symbols
follows from Theorem~\ref{thm:Manin-presentation}.
\end{proof}

Let $B^{r,N}_{x,b}$ be the transpose of the box defining $e^{r,N}_{x,b}$.
Put $m=k-j$ and
\[
 c_{r,N,k,j}=\frac{(-2\pi i)^k}{(m-1)!}r^{2-2m}N^{m-2}.
\]
By Corollary~\ref{cor:scaled-coset-comparison} and $\scrS_N^2=N^2$,
\begin{equation}\label{eq:Fourier-data-reconstruction}
 \mu_k^{(j)}(f^{r,N}_{x,b})
 =c_{r,N,k,j}N^2\pi_{\cC_k}H_k^{\mathrm{alg}}
       \bigl(z_{\Eis}(k,j)(B^{r,N}_{x,b})\bigr).
\end{equation}
The values with $j\ne2$ and $\langle x,b\rangle\not\equiv0\pmod N$
include all the endpoint data of Corollary~\ref{cor:endpoint-reconstruction}.

\section{Reconstruction from the Beilinson--Kato distribution}\label{sec:BK}

We now obtain the values in \eqref{eq:Fourier-data-reconstruction} from
the Beilinson--Kato distribution.  The local reciprocity formula applies
only to certain boxes.  We first remove its smoothing factors at a fixed
level, then use the transformation law of the Eisenstein series to pass
to the boxes required by the uniqueness theorem.

\subsection{Siegel units and the Beilinson--Kato distribution}

Let $K$ be the union of the fields of modular functions over all finite
levels, and put $\widetilde K_2^M(K)=K_2^M(K)\otimes\Q$.
On test functions vanishing in a neighbourhood of zero, let $g$ be
the Siegel distribution of \cite[Theorem~1.8]{Colmez}.  Its value on a
rational box avoiding zero is
\[
 g([a+r\Z^2])=g_{r^{-1}a},
\]
where the right side is the rational class of the Siegel unit.  These are
the values used in \cite[Proposition~1.1.2]{BPPS}.
For two such test functions, the ordered product in the row convention is
\[
 \boldsymbol\mu_2^{\mathrm{row}}(\varphi_1\otimes\varphi_2)
 =\{g(\varphi_1),g(\varphi_2)\}.
\]
By \cite[Definition~1.1.3 and Theorem~1.1.4]{BPPS}, this is the value of
their modular symbol $\xi_{\mathrm{BPPS}}$ on $[0]-[\infty]$, restricted
to these tests.  Transposition $t$ identifies row tests with column tests;
on the corresponding restricted test space,
\begin{equation}\label{eq:partII-orientation}
 (-t_*\xi_{\mathrm{BPPS}})(D_\infty)
 =t_*\boldsymbol\mu_2^{\mathrm{row}}.
\end{equation}
The two row tests specify the two factors in the displayed Milnor symbol,
in order.  We retain this presentation for the arithmetic construction.
The algebraic $q$-expansions supply the Siegel distribution to which
Wang applies the $p$-adic Kummer map, and the two rows specify the order
of the cup product.  Thus the input to Theorem~\ref{thm:intro-comparison}
includes the defining Siegel factors.  Independence of the subsequent
choices of local box and smoothing parameters is proved in
Proposition~\ref{prop:arithmetic-box-value}.

\subsection{Local reciprocity and finite-level smoothing}

Fix a prime $p$ and an embedding $\overline\Q\hookrightarrow\overline\Q_p$,
with $\overline\Q\subset\C$ as in \cite[\S1.1]{Wang}.  Wang's local
formula applies to rational test functions whose supports have $p$-adic
closure contained in $\G(\Z_p)$.  Such a test is a finite linear
combination of characteristic functions of rational boxes with this
property; these are the tests in \cite[Theorem~1.4, p.~683]{Wang}.

The ordered Kummer cup of the two Siegel distributions is Wang's
$z_{\mathrm{Kato}}$ \cite[\S2.3.1, p.~705]{Wang}.  Apply the two
smoothing operators of \cite[\S2.3.2, pp.~705--707]{Wang}, followed by
the moment
\begin{equation}\label{eq:Wang-moment}
 x\longmapsto(x_{11}e_1+x_{12}e_2)^{k-2}(\det(x)t)^{-j}.
\end{equation}
Here $x$ is the $p$-component, $e_1,e_2$ are the standard basis, and
$t$ is the compatible system of $p$-power roots of unity fixed in
\cite[\S2.3.3, p.~707]{Wang}.  Denote by
$R^W_{p;c,d;k,j}(H)$ the value on a test $H$ of Wang's local exponential
of this class.  The exponential is the construction of
\cite[\S5.1, pp.~739--742]{Wang}.

For $m=k-j$, Wang's formula reads
\begin{equation}\label{eq:Wang-reciprocity}
 R^W_{p;c,d;k,j}(H)
 =\frac1{(j-1)!}
   ({}_cz'_{\Eis}(m)\otimes{}_dz_{\Eis}(j))(H).
\end{equation}
This is \cite[Theorem~1.4 and Proposition~2.14]{Wang}; the right side is
viewed through its algebraic $q$-expansion.  Wang's product convention
omits the factor $(-1)^j$ in Definition~\ref{def:Colmez-product}.

Choose a finite decomposition of $H$ into row rectangles with a common
period, and let $L$ annihilate all their torsion parameters.  Choose
integers $c,d>1$ with $c,d\equiv1\pmod{6Lp}$, and put
\begin{equation}\label{eq:finite-smoothing-scalar}
 s_{c,d;k,j}=(c^2-c^{2-m})(d^2-d^j).
\end{equation}
These congruences fix the torsion parameters in each rectangle, so
the two smoothing operators act there by scalars.
Wang's formulas (2.5), p.~697, give
\[
 {}_cF^{(m)}=(c^2-c^{2-m})F^{(m)},\qquad
 {}_dE^{(j)}=(d^2-d^j)E^{(j)}.
\]
When $m=2$, the first parameter is nonzero on every such rectangle.
When $j\ne2$, both scalar factors are nonzero.

\begin{proposition}\label{prop:local-unsmoothing}
For $k\ge2$, $1\le j\le k-1$, $j\ne2$, and a test $H$ as above,
\begin{equation}\label{eq:local-unsmoothed-value}
 F_H:=\frac{(-1)^j}{s_{c,d;k,j}}R^W_{p;c,d;k,j}(H)
     =z_{\Eis}(k,j)(H).
\end{equation}
In particular, $F_H$ is independent of the common decomposition and of
$c,d$ satisfying the stated congruences.
\end{proposition}

\begin{proof}
Substitute the two scalar identities in \eqref{eq:Wang-reciprocity}.
The right side of \eqref{eq:local-unsmoothed-value} is the algebraic
modular form specified in Definition~\ref{def:Colmez-product}; thus the
identity of $q$-expansions specifies that form and its normalization.
\end{proof}

The theorem applies also at $p=2$.  The auxiliary level used in Wang's
calculation may be refined to satisfy $v_p(M)\ge v_p(2p)$
\cite[\S1.2.4 and Proposition~5.5]{Wang}.  This refinement does not
change \eqref{eq:local-unsmoothed-value}.

\subsection{Rational changes of box}

The Fourier values used in Theorem~\ref{thm:Fourier-Manin-uniqueness}
need not correspond to boxes in $\G(\Z_p)$ for a fixed prime $p$.
We shall choose $p$ and a rational matrix separately for each box.
Once reciprocity has given an Eisenstein form, its usual transformation
law determines the value on the original box.

For $g=\left(\begin{smallmatrix}a&b\\c&d\end{smallmatrix}\right)
\in\G^+(\Q)$ and a weight-$h$ form $f$, put
\[
 (f\star_hg)(\tau)=(c\tau+d)^{-h}
 f\left(\frac{a\tau+b}{c\tau+d}\right).
\]
This is the algebra action of \cite[(1.2), p.~679]{Wang}.
It is used here only for the intermediate Eisenstein forms.
The one-factor identities of \cite[Theorem~2.13, pp.~695--696]{Wang}
give, for rational row tests $U$ and $\delta=\det g>0$,
\[
 z_{\Eis}(j)(Ug)=z_{\Eis}(j)(U)\star_jg,\qquad
 z'_{\Eis}(m)(Ug)=\delta^{m-1}z'_{\Eis}(m)(U)\star_mg.
\]
Here $j\ne2$; the primed weight-two identity is used on the nonzero
parameters occurring below.  Multiplication of forms therefore gives,
with $Z=z_{\Eis}(k,j)$,
\begin{equation}\label{eq:Eisenstein-box-transport}
 Z(Bg)=\delta^{m-1}(Z(B)\star_kg),\qquad
 Z(B)=\delta^{1-m}(Z(Bg)\star_kg^{-1}).
\end{equation}
These identities follow first for row rectangles and then by linearity.

\begin{lemma}\label{lem:box-unit-normalization}
Let $B=(r/N)A+r\Mat(\Z)$, where $r\in\Q_{>0}$,
$A\in\Mat(\Z)$, and $\det A\not\equiv0\pmod N$.
There are a prime $p\mid N$ and $g\in\G^+(\Q)$ such that the
$p$-adic closure of $Bg$ lies in $\G(\Z_p)$.
\end{lemma}

\begin{proof}
Choose $p\mid N$ with $v_p(\det A)<v_p(N)$ and put
$g=(N/r)A^{-1}$.  Then
\[
 Bg=I+N\Mat(\Z)A^{-1},\qquad
 NA^{-1}=N\operatorname{adj}(A)/\det A\in p\Mat(\Z_p).
\]
Thus the closure lies in $I+p\Mat(\Z_p)$.  If $\det g<0$, replace
$g$ by $g\operatorname{diag}(-1,1)$; the closure remains in
$\G(\Z_p)$.  The rational lattice coset $Bg$ is a finite disjoint
union of scalar-lattice boxes, since its lattice contains
$s\Mat(\Z)$ with finite index for some $s\in\Q_{>0}$.
\end{proof}

\begin{proposition}\label{prop:arithmetic-box-value}
Let $k\ge2$, $1\le j\le k-1$, $j\ne2$, and $m=k-j$.
For a box $B$ as in Lemma~\ref{lem:box-unit-normalization}, choose
$p,g$ as in that lemma and compute $F_{Bg}$ by \eqref{eq:local-unsmoothed-value}.
Then
\begin{equation}\label{eq:arithmetic-box-value}
 Z_B:=(\det g)^{1-m}(F_{Bg}\star_kg^{-1})
       =z_{\Eis}(k,j)(B).
\end{equation}
The expression is independent of $p,g,c,d$ and of the decomposition of
$Bg$.  It respects finite refinements of these boxes.
\end{proposition}

\begin{proof}
Proposition~\ref{prop:local-unsmoothing} identifies $F_{Bg}$ with
$Z(Bg)$.  Apply \eqref{eq:Eisenstein-box-transport}.  Independence and
refinement compatibility follow from the resulting equality with $Z(B)$.
\end{proof}

\subsection{Reconstruction}

\begin{proof}[Proof of Theorem~\ref{thm:intro-comparison}]
For $\langle x,b\rangle\not\equiv0\pmod N$, apply
Proposition~\ref{prop:arithmetic-box-value} to the row box
$B^{r,N}_{x,b}$.  Equations \eqref{eq:arithmetic-box-value} and
\eqref{eq:Fourier-data-reconstruction} give
\[
 \mu_k^{(j)}(f^{r,N}_{x,b})
 =c_{r,N,k,j}N^2\pi_{\cC_k}H_k^{\mathrm{alg}}(Z_{B^{r,N}_{x,b}})
 \qquad(j\ne2).
\]
Thus the analytically constructed $\mu$ has the prescribed values.
In each weight these data include the coefficients required by
Theorem~\ref{thm:Fourier-Manin-uniqueness}, since their indices all
satisfy $j\ne2$.  This gives uniqueness of $\mu$ and hence of $\Phi$
by evaluation at $D_\infty$.
\end{proof}

This reconstruction uses the local exponential in
\eqref{eq:Wang-reciprocity}.  Its identification with the geometric
Bloch--Kato exponential is not required; compare
\cite[\S1.2.5, p.~684]{Wang}.

\begin{remark}
At level $N=6$, take $x=(2,0)$ and $b=(0,1)$.  Then
$\langle x,b\rangle=2$ is nonzero but not a unit modulo $6$.
The corresponding row box is
$B=\operatorname{diag}(1/3,1/6)+\Mat(\Z)$.
For $p=3$ and $g=\operatorname{diag}(3,6)$, its image is
\[
 Bg=I+\Mat(\Z)\operatorname{diag}(3,6)
 \subset I+3\Mat(\Z_3).
\]
In weight three, reciprocity determines the coefficient with $j=1$;
the coefficient with $j=2$ is recovered from
$\mu_3^{(2)}(U)=-\mu_3^{(1)}(US)$.
\end{remark}

\section{Relations among Eisenstein products}
\label{sec:level-one}

Evaluating the cyclic identity on the integral lattice gives a
polynomial identity among level-one Eisenstein series.  We determine
the space of relations obtained by taking its coefficients.

For \(m\ge1\), put \(\mathbb E_m=\phi_{m,\Z^2}(0,\tau)\) and
\(\Pi_n(X,Y)=X^n+Y^n+(-X-Y)^n\).  The all-integral specialization of
\eqref{eq:all-integral-square}, in homogeneous degree \(n>0\), is
\begin{align}\label{eq:journal-full-triangle}
 &\sum_{i+j=n}\Pi_i\Pi_j\mathbb E_{i+1}\mathbb E_{j+1}\notag\\
 &\qquad=\left((n+1)\Pi_n+
 \frac{\Pi_{n+1}\Pi_2}{XY(X+Y)}\right)\mathbb E_{n+2}.
\end{align}
For even \(n\) the quotient is a polynomial; for odd \(n\) both sides
vanish.

\begin{lemma}[Coefficient form of the level-one identity]
\label{lem:level-one-HST}
For $k=2h\ge8$, the coefficient relations of
\eqref{eq:journal-full-triangle}, with $n=k-2$, generate all
linear relations among \(\mathbb E_{2h}\) and the unordered products
\(\mathbb E_{2r}\mathbb E_{2h-2r}\),
\(2\le r\le\lfloor h/2\rfloor\).
\end{lemma}

\begin{proof}
For odd \(m\),
\begin{equation}\label{eq:Pi-binomial}
 \Pi_m(X,Y)
 =-\sum_{a=1}^{m-1}\binom maX^aY^{m-a}.
\end{equation}
Put $q=\lfloor(k-2)/6\rfloor$.  For $2\le a\le q+1$, define the
linear functional on homogeneous polynomials of degree $k-2$ by
\[
 \lambda_a(P)=[U^{2a-2}V^{2a-2}W^{k-4a+2}]P(U-V,V-W).
\]
Write $(A,B,C)=(2a-2,2a-2,k-4a+2)$.  These are positive even integers,
and $C\ge A=B$.  For positive odd $i,j$ with $i+j=k-2$, the binomial
theorem gives
\[
 \lambda_a(\Pi_i\Pi_j)
 =-\sum_{(d,e,f)}\binom{i}{d}\binom{j}{e},
\]
where the sum runs over the six permutations of the three positions
in $(A,B,C)$, with repetitions, and $\binom{i}{d}=0$ for $d>i$.
Indeed, the terms $(U-V)^{i+j}$, $(V-W)^{i+j}$ and $(W-U)^{i+j}$
have zero coefficient because $A,B,C>0$.  Each of the six cross terms
has the negative sign displayed above.

Apply $\lambda_a$ to \eqref{eq:journal-full-triangle}.  The coefficient
of $\mathbb E_{2b}\mathbb E_{k-2b}$ vanishes for $b<a$, since
$2b-1<A$.  It is strictly negative for $b=a$: every summand is
nonpositive and, in particular,
\[
 \binom{2a-1}{A}\binom{k-2a-1}{C}>0.
\]
The resulting $q$ relations are therefore independent, with successive
pivots on the products indexed by $a=2,\ldots,q+1$.

By \cite[Corollary~2]{HiroseSatoTasaka}, $\mathbb E_k$ together with
the products indexed by $q+2,\ldots,\lfloor k/4\rfloor$ is a basis
of $M_k(\SL(\Z))$; the normalization change is
$\mathbb E_m=(2\pi i)^m E_m^{\mathrm{HST}}$ for even $m\ge4$.
The relation space therefore has dimension $q$.  The $q$ independent
relations obtained by the $\lambda_a$ belong to the span of the
coefficient relations and exhaust this space.
\end{proof}

\begin{theorem}[Level-one product relations]
\label{thm:journal-level-one-relations}
Let \(k=2h\ge8\), and let
\[
 m_{2h}:\bigoplus_{r=2}^{\lfloor h/2\rfloor}\C T_r
 \longrightarrow M_{2h}(\SL(\Z)),\qquad
 T_r\longmapsto\mathbb E_{2r}\mathbb E_{2h-2r}.
\]
Then \(\ker(m_{2h})\) consists of the linear combinations of the coefficient
relations obtained from \eqref{eq:journal-full-triangle} with \(n=2h-2\)
whose total \(\mathbb E_{2h}\)-coefficient is zero.  In particular,
\[
 \dim\ker(m_{2h})=\left\lfloor\frac{h-4}{3}\right\rfloor
 =\left\lfloor\frac{k-8}{6}\right\rfloor.
\]
\end{theorem}

\begin{proof}
Adjoin \(T_0\), mapped to \(\mathbb E_{2h}\), to the source of
\(m_{2h}\).  By Lemma~\ref{lem:level-one-HST}, the coefficient relations
of \eqref{eq:journal-full-triangle} generate the kernel of this enlarged
map, of dimension $q=\lfloor(k-2)/6\rfloor$.  The kernel of $m_{2h}$
is its subspace with zero $T_0$-coefficient.

The coefficient of \(\mathbb E_{2h}\) in
\eqref{eq:journal-full-triangle} is
\[
 C_h(X,Y)=(2h-1)\Pi_{2h-2}
 +\frac{\Pi_{2h-1}\Pi_2}{XY(X+Y)}.
\]
It is nonzero, since
\[
 C_h(1,1)=4h+4+(2h-7)2^{2h-2}>0.
\]
Thus the $T_0$-coefficient is a nonzero linear functional on the enlarged
kernel.  Its kernel has dimension
$q-1=\lfloor(k-8)/6\rfloor$.
\end{proof}

At higher level, multiplying identities requires some care:
holomorphic projection is linear, but not multiplicative.  We therefore
retain the derivative terms which make each projected quadratic identity
an identity of nearly holomorphic forms before multiplying it by another
form.  Fix $N$, and put
\[
 V_{N,m}=\operatorname{span}_{\C}
 \{E_{m,a}:a\in N^{-1}\Z^2/\Z^2\},\qquad
 V_N=\bigoplus_{m\ge1}V_{N,m},
\]
with the completed interpretation in weight two.  We need a polynomial
ring in which an Eisenstein series and its first raised derivative are
represented by independent variables.  Linear relations among the
Eisenstein series themselves are to remain imposed.  Choose a homogeneous
basis $\mathcal B_N$ of $V_N$ and form
\[
 \mathcal A_N=\C[x_f,y_f:f\in\mathcal B_N],\qquad
 \operatorname{wt}(x_f)=m,\quad\operatorname{wt}(y_f)=m+2
 \quad(f\in V_{N,m}).
\]
Here $x_f,y_f$ represent $f$ and
$D_mf=\delta_mf/(2\pi i)$, respectively.  For arbitrary $f\in V_N$,
extend $x_f$ and $y_f$ linearly.  Thus linear relations among Eisenstein
series are imposed, but identities involving their products or derivatives
are not.

For $f\in V_{N,2}$, write
$f=f^{\mathrm{hol}}+\epsilon(f)/(4\pi\operatorname{Im}\tau)$;
put $\epsilon(f)=0$ in other weights.  For homogeneous $f,g$ of weights
$m,n$, define
\begin{equation}\label{eq:quadratic-formal-lift}
 \mathcal L_N(fg)=x_fx_g+
       \frac{\epsilon(f)}{n}y_g+
       \frac{\epsilon(g)}{m}y_f,
\end{equation}
and extend linearly to formal quadratic expressions.  Substitution
$x_f=f$, $y_f=D_mf$ sends $\mathcal L_N(fg)$ to
$H_{m+n}^{\mathrm{alg}}(fg)$: the nonholomorphic terms cancel, and the
correction is a sum of raised forms.

Accordingly, a homogeneous relation modulo Eisenstein series,
\[
 H_k^{\mathrm{alg}}(Q)=e,\qquad e\in V_{N,k},
\]
is represented by the polynomial $\mathcal L_N(Q)-x_e$ in
$\mathcal A_N$.  The homogeneous coefficients of the cyclic $R$-identity
give such relations by Proposition~\ref{prop:Manin2}, including its
regularized zero-parameter cases.

Let $\mathcal I_N\subset\mathcal A_N$ be the ideal generated by the
cyclic relations.  Since each generator has weight at least two, the
weight-two part of $\mathcal I_N$ is the linear span of its weight-two
generators.  To determine the equations on the projective charts
$x_f\ne0$, $f\in V_{N,1}$, one must allow multiplication by powers of
the chart coordinate.  This leads to the following saturation result.

There is a geometric precedent for this distinction.
Borisov, Gunnells, and Popescu
\cite[Theorem~4.4 and Corollary~4.5]{BGP} give scheme-theoretic equations
for the embedding of $X_1(p)$, $p\ge5$ prime, defined by weight-one Eisenstein series,
while distinguishing this conclusion from generation of its homogeneous
ideal \cite[Remark~4.6(3)]{BGP}.  In the full-level weight-one embedding
for which Khuri-Makdisi obtains equations of degrees two and three
\cite[end of \S5]{KhuriMakdisi}, we can specify a defining family of quadrics.

Write $G=(\Z/p\Z)^2$, $\zeta=e^{2\pi i/p}$, and
$(Fh)_a=\sum_{b\in G}\zeta^{a_1b_2-a_2b_1}h_b$.
Use the weight-one coordinates $s_a=-E_{1,a/p}/(2\pi i)$ for $a\ne0$,
with $s_0=0$.  Their linear relations are $s_{-a}=-s_a$ and $Fs=-ps$.

\begin{theorem}[Equations for the full prime-level curve]\label{thm:prime-weight-two-saturation}
Let $p$ be an odd prime.  Put $v_a=s_a^2$ and $u=2v+Fv$.  The image
of $X(p)$ in $\mathbf P(V_{p,1}^{\vee})$ over $\C$ is cut out
scheme-theoretically by the quadrics
\begin{equation}\label{eq:prime-defining-quadrics}
 (p^2-4)\bigl(s_as_b-s_{a+b}(s_a+s_b)\bigr)
 =u_a+u_b+u_{a+b}
 \qquad(a,b,a+b\ne0).
\end{equation}

In particular, if $Q$ is a quadratic expression in $V_{p,1}$
and $Q=e$ as modular forms, with $e\in V_{p,2}$, then there is an integer
$r\ge0$ such that
\[
 x_{f_1}\cdots x_{f_r}\bigl(\mathcal L_p(Q)-x_e\bigr)
       \in\mathcal I_p
 \qquad(f_1,\ldots,f_r\in V_{p,1}).
\]
Only the weight-two generators of $\mathcal I_p$ are needed.
\end{theorem}

Substituting $t=-u/(p^2-4)$ from \eqref{eq:prime-elimination} in the
weight-two cyclic Manin relation \eqref{eq:prime-cyclic} gives
\eqref{eq:prime-defining-quadrics}.  Appendix~\ref{app:prime-weight-two}
proves the scheme equality and deduces the ideal membership by saturation.  In weight two,
$\mathcal L_p(fg)=x_fx_g$ for $f,g\in V_{p,1}$, and $e$ is holomorphic.
This suggests the following statement in arbitrary level and weight.

\begin{conjecture}[Relations after saturation]
\label{conj:R-relation-ideal}
Let $N\ge3$.  For every homogeneous quadratic relation
$H_k^{\mathrm{alg}}(Q)=e$ with $e\in V_{N,k}$, there exists $r\ge0$ such that
\[
 x_{f_1}\cdots x_{f_r}\bigl(\mathcal L_N(Q)-x_e\bigr)
       \in\mathcal I_N
 \qquad\text{for all }f_1,\ldots,f_r\in V_{N,1}.
\]
\end{conjecture}

Equivalently, each relation belongs to $\mathcal I_N$ after localizing
at any weight-one basis variable.  Indeed, the displayed condition
implies local membership.  Conversely, if powers of each of the finitely
many weight-one basis variables carry the relation into $\mathcal I_N$,
every monomial of sufficiently large degree in those variables does so.
Thus the conjecture asserts ideal membership away from the locus
where all weight-one coordinates vanish.  The restriction $N\ge3$ is
necessary for this formulation:
at levels one and two the weight-one Eisenstein space is zero.

In higher weights the assertion concerns the corrected expressions
$\mathcal L_N(Q)$, not products of already projected relations.
Neither ordinary ideal generation nor a uniform bound on $r$ is asserted.

\appendix
\section{Weight-two relations at prime level}
\label{app:prime-weight-two}

We prove Theorem~\ref{thm:prime-weight-two-saturation}.  All rings and
schemes in this appendix are over $\C$, except where a coefficient
ring is explicitly specified.  We first eliminate the weight-two
coordinates from the cyclic equations.  To a solution of the resulting
equations we associate a curve, following \cite[\S\S3--4]{BGP}.
Away from the cusps this curve is elliptic, and its torsion points
recover the given Eisenstein coordinates.  We then check first-order
deformations, both away from the cusps and at the cusps, to obtain
equality of schemes rather than just equality of their complex points.
At full level the argument must also establish connectedness and the
prescribed Weil pairing; neither follows from restriction to a cyclic
subgroup.

\subsection{The equations in Eisenstein coordinates}

The weight-one coordinates will be the Eisenstein series already used
in the paper.  For the auxiliary elliptic curve it is convenient to
use their Laurent normalization, as in \cite{BGP,KhuriMakdisi}.
Write $G=(\Z/p\Z)^2$, $M=p^2$, and $\zeta=e^{2\pi i/p}$.
For $a\ne0$ put
\[
 s_a=-\frac{E_{1,a/p}}{2\pi i},\qquad
 w_a=\frac{E_{2,a/p}-E_{2,0}}{(2\pi i)^2},\qquad
 t_a=-\frac1M\sum_{b\in G}\zeta^{b_2a_1-b_1a_2}w_b,
\]
where $s_0=t_0=w_0=0$.  Here $w_a$ is the normalized Weierstrass
coordinate; the $t_a$ are another family spanning the
holomorphic weight-two Eisenstein space.  Define the finite Fourier
transform by
\[
 (Fh)_a=\sum_{b\in G}\zeta^{b_2a_1-b_1a_2}h_b.
\]
The identities of \cite[(3.28), (3.29) and Proposition~4.3]{KhuriMakdisi}
give
\begin{equation}\label{eq:prime-linear}
 s_{-a}=-s_a,\quad Fs=-ps,\qquad
 t_{-a}=t_a,\quad \sum_a t_a=0,\qquad w=-Ft,
\end{equation}
and
\begin{equation}\label{eq:prime-cyclic}
 C(a,b):=s_as_b-s_{a+b}(s_a+s_b)+t_a+t_b+t_{a+b}=0
 \quad(a,b,a+b\ne0).
\end{equation}
Moreover, $s_a^2=w_a+2t_a$.  The odd Fourier eigenspace in
\eqref{eq:prime-linear} has dimension $(M-1)/4$, and the even $t$-space
has dimension $(M-3)/2$.  These are the dimensions of the corresponding
Eisenstein spaces, so these are their full linear relations.
We henceforth regard $s_a,t_a$ as coordinates subject to these linear
relations, not as functions already restricted to the modular curve.
Equation~\eqref{eq:prime-cyclic} is precisely the weight-two cyclic
relation in these coordinates.  The zero-parameter cases give the
square relations $s_a^2=w_a+2t_a$.

In fact the square relations follow from the nonzero cyclic equations.
Put $v_a=s_a^2$.  Since $F$ is skew-transpose on odd vectors,
$s^tFs=0$; hence $\sum_a v_a=0$.  Fourier transformation of convolution,
using $Fs=-ps$ and $F^2=M$, gives $s*s=Fv$.  Sum
\eqref{eq:prime-cyclic} over $b\ne0,-a$.  The result is
\begin{equation}\label{eq:prime-elimination}
 t=-\frac{2v+Fv}{M-4},\qquad
 w=-Ft=\frac{Mv+2Fv}{M-4}=v-2t.
\end{equation}
This calculation takes place in the polynomial quotient by the cyclic
equations.  It eliminates all holomorphic weight-two coordinates.

Let $S$ be the projective scheme in the weight-one Eisenstein space
defined by substituting \eqref{eq:prime-elimination} into
\eqref{eq:prime-cyclic}.  The actual Eisenstein series give a closed
embedding $X(p)\hookrightarrow S$.  To see that this is an embedding,
Khuri-Makdisi's generation theorem
\cite[Theorems~3.5 and 5.1]{KhuriMakdisi} gives
$\operatorname{Sym}^2 V_{p,1}\twoheadrightarrow M_2(\Gamma(p))$.
The weight-two line bundle has degree $p(p^2-1)/12$, which exceeds
$2g(X(p))+1$ by $(p^2-1)/2-3>0$.  It is very ample, and its map factors
through the quadratic Veronese of the weight-one map.  The latter is
therefore a closed embedding.  It remains to prove that this embedding
is an equality, including the scheme structure.  For $p=3$ this is immediate:
$X(3)$ has genus zero and its weight-one line bundle has degree one.
Its two-dimensional Eisenstein space is the full space of sections,
so its image is the whole ambient projective line.  Consequently every
defining quadric of $S$ is zero and $S=X(3)$.  The elimination above
already proves the theorem in this case.  For the rest of the proof
assume $p\ge5$.

\begin{lemma}\label{lem:prime-line-restriction}
At every point of $S$, restriction of the $s_a$ to any nonzero cyclic
subgroup of $G$ is nonzero.
\end{lemma}
\begin{proof}
By a symplectic change of indices it suffices to consider
$H=\{(0,k):k\in\Z/p\Z\}$.  Suppose $s|_H=0$.
Summing the cyclic equations within $H$ shows that its $t$-coordinates
are constant; an individual equation then shows that they vanish.
For $r\ne0$ put $u_k=s_{(r,k)}$ and $v_k=t_{(r,k)}$.
The mixed equations give $u_i u_j=v_i+v_j$ for $i\ne j$.
Thus, for four distinct indices,
$(u_i-u_k)(u_j-u_l)=0$.  Since $p\ge5$, all entries of a row are equal
with at most one exception.  The equation $Fs=-ps$ on $H$, together
with oddness, makes every row sum zero.  Consequently
$u_k=c_r(1-p\delta_{k,k_r})$ for some $c_r,k_r$.
For $R\ne0$ the same Fourier equation becomes
\[
 \sum_{r\ne0}c_r\zeta^{k_rR-rT}
       =c_R(1-p\delta_{T,k_R}).
\]
The coefficients of $\zeta^{-RT}$ on the two sides are
$c_R\zeta^{k_RR}$ and its negative.  Hence $c_R=0$ for every $R$,
contradicting the fact that $s$ is a projective point.
\end{proof}

The restriction to a line satisfies the equations for $X_1(p)$ in
\cite[Theorem~4.4]{BGP}.  By the lemma it is a projective point, so
their theorem applies to every such restriction.  To recover a point
of $X(p)$, however, we must construct a single elliptic curve carrying
the full level structure.

\subsection{The auxiliary curve}

We next construct an elliptic curve from a solution of
\eqref{eq:prime-cyclic}.  We give the algebraic details also over a
coefficient ring, since they will be needed for first-order deformations.
In this subsection only oddness of $s$, evenness of $t$, and the cyclic
equations are used; the Fourier relations are not required.

For $a\ne0$ there are unique formal Laurent series
\begin{equation}\label{eq:prime-ode}
 \widehat r_a(z)=z^{-1}+s_a+t_a z+u_a z^2+\cdots,
 \qquad
 \widehat r_a'=-\frac{\sum_{k\ne0,a}\widehat r_k\widehat r_{a-k}
                              +2s_a\widehat r_a}{M-2}.
\end{equation}
The cyclic equations give consistency at the first two orders.
For the coefficient of $z^n$ the recursion has matrix
\[
 \left(n-\frac2{M-2}\right)I+\frac2{M-2}J,
\]
where $J$ is the all-ones matrix on the $M-1$ indices.
Its eigenvalues are $n-2/(M-2)$ and $n+2$, both nonzero for $n\ge2$.
This proves existence and uniqueness over the universal coefficient
ring, and hence after any specialization.  Uniqueness also gives
$\widehat r_a(-z)=-\widehat r_{-a}(z)$, and
\begin{equation}\label{eq:prime-u}
 (M-3)u_a=\sum_{k\ne0,a}s_{k-a}t_k-s_at_a.
\end{equation}

Introduce homogeneous coordinates $r_0,r_a$ for $a\ne0$.  The curve
will be defined by the quadrics
\begin{align}
 R_{a,b;c,d}&=r_ar_b-r_cr_d-r_0r_{a+b}(s_a+s_b-s_c-s_d),
       &&a+b=c+d\ne0,\label{eq:prime-aux-nonzero}\\
 R^0_{a,b}&=r_ar_{-a}-r_br_{-b}
                    +r_0^2(s_a^2-2t_a-s_b^2+2t_b).
       \label{eq:prime-aux-zero}
\end{align}
All subscripts other than $0$ in these formulas are nonzero.
These are the auxiliary equations of \cite[\S3]{BGP}, now indexed by
$G$ rather than by a cyclic group.

\begin{lemma}\label{lem:prime-auxiliary}
The substitution $r_0=1$, $r_a=\widehat r_a(z)$ annihilates
\eqref{eq:prime-aux-nonzero}--\eqref{eq:prime-aux-zero}.
Over any coefficient ring satisfying the cyclic equations, the degree
$n$ part of their homogeneous quotient is free of rank $nM$ for $n\ge1$.
Over $\C$ their projective zero scheme is a reduced curve, of degree $M$
and arithmetic genus one.  Its intersection with $r_0=0$ consists of the
$M$ distinct smooth points
\[
 P_\chi=[0:(\chi(a))_{a\ne0}],\qquad \chi\in G^*.
\]
Every component contains one of these points.
\end{lemma}
\begin{proof}
We first verify the Laurent substitution.  The constant terms of the
quadrics vanish by \eqref{eq:prime-cyclic}.  For the next coefficient
the required identity, with $d=a+b\ne0$, is
\[
 s_at_b+s_bt_a-(s_a+s_b)t_d+u_a+u_b+2u_d=0.
\]
To verify this over the coefficient ring, multiply the left side by
$M-3$ and express it in terms of the cyclic equations:
\begin{align*}
 &\sum_{k\ne0,-a,b}\bigl\{
 (s_{k-b}+s_b)C(-k-a,k)+(s_a-s_{k+a})C(k-b,-k)\\
 &\hspace{45mm}-(s_a+s_b)C(b-k,k+a)\bigr\}\\
 &\hspace{15mm}-(s_a+s_b+2s_d)C(a,b).
\end{align*}
Indeed, substitute \eqref{eq:prime-u}, expand, and collect the terms
in $t_k,t_{k+a},t_{k-b}$.  Reindexing these sums and using oddness
cancels their complete sums; the excluded indices give the last
displayed multiple of $C(a,b)$.

The differential equation preserves the auxiliary quadrics.  More
explicitly, on $r_0=1$, for $e=a+b=c+d\ne0$ it gives
\begin{align*}
 (M-2)R_{a,b;c,d}'={}&
 -\sum_{k\ne0,a,c,e}r_kR_{a-k,b;c-k,d}
 -\sum_{k\ne0,b,d,e}r_kR_{b-k,a;d-k,c}\\
 &-r_e(R^0_{a,c}+R^0_{b,d})
 -(s_a+s_b+s_c+s_d)R_{a,b;c,d},\\
 (M-2)(R^0_{a,b})'={}&
 -\sum_{k\ne0,a,b}\bigl(r_kR_{-a,a-k;-b,b-k}
                    +r_{-k}R_{a,k-a;b,k-b}\bigr).
\end{align*}
To check these identities, substitute \eqref{eq:prime-ode} into each
derivative and group terms with common factor $r_k$.
Before imposing the cyclic equations, the first identity has the
additional term $-2r_e(C(a,b)-C(c,d))$; the second has no remainder.
The coefficient of $z$ in $R^0$ vanishes by parity.  If the first
nonzero coefficient of any of the substituted quadrics has degree $m\ge2$,
these differential identities give a linear system with diagonal
$m(M-2)$ and total absolute right-hand coefficient at most $2M-6$.
Identically zero quadrics are omitted.  Since $m(M-2)>2M-6$, the matrix
is invertible.  Induction proves the asserted vanishing.

We next compute the Hilbert function.  Grade the polynomial ring also by
$G$, giving $r_a$ degree $a$ and $r_0$ degree zero.  Modulo $r_0$, the
auxiliary equations allow one pair of factors to be replaced by another
having the same sum of indices.  We claim that these replacements
identify any two monomials with the same ordinary degree and the same
$G$-degree.  A pair $r_ar_b$ may acquire a specified
factor $r_c$ unless $a+b=c$.  If two monomials of degree at least three
have no common factor and no pair replacement can give them one, the
sum of any pair of indices in the first must equal every index occurring
in the second, and conversely.  Each monomial is then
a pure power, say of $r_a$ and $r_b$, with $2a=b$ and $2b=a$.
This would give $3a=0$, impossible for $p\ge5$.
Once a factor is shared, induction reduces the degree.  Degree two is
a defining relation.  Multiplying inductively by $r_0$ now gives at
most $n$ generators in each $G$-degree.

To prove independence, in a nonzero $G$-degree choose monomials
of pole orders $1,\ldots,n$, and in degree zero choose pole orders
$2,\ldots,n$ and $r_0^n$.  Such monomials exist: a nonzero index can
always be split into two nonzero indices.  Multiply by powers of $r_0$
to make ordinary degree $n$.  Their Laurent images have distinct
leading powers with coefficient one, so they are independent over the
coefficient ring.  The preceding spanning argument can use these same
monomials.  They are therefore a free basis.

The Laurent series also describe the components of the projective
scheme.  To keep different homogeneous degrees separate, introduce an
indeterminate $T$ and use all substitutions
\[
 r_0\longmapsto T,\qquad
 r_a\longmapsto T\chi(a)\widehat r_a(z),\qquad\chi\in G^*.
\]
Character orthogonality separates the $G$-degrees, and $T$ separates
the ordinary degrees.  The independence just proved therefore embeds
the graded quotient over $\C$ into
$\prod_{\chi\in G^*}\C((z))[T]$.  In particular, it is reduced.
Each substitution has prime kernel.  Since their intersection is zero,
every minimal prime is one of these kernels.  Its projective image is
not a point: the ratio $r_0/r_a$ maps to
$\chi(a)^{-1}z+O(z^2)$ and is nonconstant.  On the other hand, the
Hilbert function $Mn$ shows that the whole projective scheme has
dimension one.  Thus every irreducible component is a curve; there are
no isolated components.  Multiplying the Laurent coordinates by $z$
and setting $z=0$ shows that the branch indexed by $\chi$ meets
$P_\chi$.  The Hilbert polynomial is $Mn$, so the curve has degree $M$
and arithmetic genus one.
Multiplication by $r_0$ is injective, and its hyperplane section has
length $M$.  The $M$ displayed points already exhaust this length;
the section is reduced.  At each of them the local maximal ideal is
generated by the hyperplane parameter.  The curve is therefore smooth
there.  The free-basis argument remains valid over the dual numbers.
\end{proof}

\subsection{Reconstruction away from the cusps}

We first show that a singular auxiliary curve can occur only at a
cusp.  We then recover an elliptic curve with full level structure from
every other solution of the cyclic equations.

Let $C$ be the auxiliary curve attached to a complex point of $S$.
Equation~\eqref{eq:prime-ode} defines a vector field on $r_0\ne0$.
At $P_\chi$ the Laurent parameter makes it $d/dz$, so it extends
regularly and is nonzero there.  If $C$ is singular, its vector field
vanishes at each singular point: a local flow preserves the finite
singular locus.  Such a point is in $r_0\ne0$; write it
$[1:(\rho_a)]$.  The auxiliary equations and the vanishing of the
vector field give
\begin{equation}\label{eq:prime-singular}
 \rho_a\rho_b=\rho_{a+b}(s_a+s_b)\qquad(a,b,a+b\ne0).
\end{equation}
Indeed, the difference between the two sides depends only on $a+b$;
summing over its decompositions makes it zero by the differential
equation.  Thus $H=\{0\}\cup\{a:\rho_a\ne0\}$ is a subgroup:
it is closed under nonzero sums, and repeated addition gives negatives.

If $H=\{0\}$, the zero-sum auxiliary equations make every $w_a$
equal.  Their total sum is zero, so they vanish.  The cyclic equations
then give $(s_a+s_b-s_{a+b})^2=0$.  Additivity on the finite group
forces $s=0$, a contradiction.

If $H=G$, the ratios $-\rho_a/\rho_{-a}$ form a character.
Twisting the $\rho_a$ by a character, using that $p$ is odd, makes
$\rho_{-a}=-\rho_a$.  For $a+b+c=0$, equation
\eqref{eq:prime-singular} gives
$\rho_a^2=(s_a+s_b)(s_a+s_c)$.  Hence
$s_as_b+s_bs_c+s_cs_a=\rho_a^2-s_a^2$ is independent of the triple.
The triples connect all nonzero indices.  Summing the cyclic equations
therefore makes all $t_a$ equal, and $\sum_a t_a=0$ makes them zero.
If some $s_a\ne0$, the equations on its cyclic subgroup give
$s_{na}=s_a/n$ for $1\le n\le p-1$.  This contradicts
$s_{-a}=-s_a$.

It remains that $H$ has order $p$.  Put it in vertical position.
For two indices outside $H$ with sum a nonzero element of $H$,
\eqref{eq:prime-singular} says that their $s$-coordinates sum to zero.
It follows that $s$ is constant on each nonzero coset of $H$.
The $w$-coordinates are all equal outside $H$, by
\eqref{eq:prime-aux-zero}.  Restrict to a transverse cyclic subgroup;
its $s$-coordinates are nonzero by Lemma~\ref{lem:prime-line-restriction}.
This is the constant-Weierstrass-coordinate case of
\cite[Lemma~4.2]{BGP}: up to scale and a nonzero reindexing, its values
are $1/2-r/p$, $1\le r\le p-1$.
The equation $Fs=-ps$ uniquely determines the vertical coordinates
from these coset constants.  They are their finite Fourier transform,
namely the cotangent constants at a cusp.  Finally,
\eqref{eq:prime-elimination} determines $t$.  Thus a singular auxiliary
curve occurs only over an actual cusp of $X(p)$.

At every other point $C$ is smooth.  The character group $G^*$ acts by
$r_a\mapsto\chi(a)r_a$ and permutes the $P_\chi$ transitively.
Since every component contains such a point, it acts transitively on
the components.  Arithmetic genus one then makes every component an
elliptic curve.  We claim that there is only one component.
Otherwise its stabilizer $K\subsetneq G^*$ acts by translations:
an origin-preserving automorphism group of an elliptic curve has no
nontrivial $p$-subgroup for $p\ge5$.  Choose a nontrivial character
$a\in G$ trivial on $K$.  On that component $r_a/r_0$ is $K$-invariant
and has simple poles exactly at the marked points, a single $K$-orbit.
It descends to the elliptic quotient with exactly one simple pole,
contrary to Riemann--Roch.  The action of $K$ is faithful, since it is
free on the marked points.  This proves connectedness.

Give $C$ origin $P_1$ and let $\omega$ be the differential dual to its
vector field.  The field is nowhere zero on this elliptic curve.
The character action is by translations and identifies $G^*$ with
$C[p]$.  Let $Q_a\in C[p]$ be determined by
\[
 e_p(Q_a,P_\chi)=\chi(a)\qquad(\chi\in G^*).
\]
For each nonzero $a$, the function $r_a/r_0$ is uniquely determined
by its character, simple poles on $C[p]$, and residue one at the
origin with respect to $\omega$.  Indeed, the difference of two such
functions has no pole at the origin, hence none at any marked point
by equivariance; it is constant, and its nontrivial character makes
that constant zero.

This also fixes the normalization of the Eisenstein coordinates.
The function $g_Q$ of \cite[Definition~4.1]{KhuriMakdisi} has character
$e_p(Q,\cdot)$ and, in the parameter $z=\int\omega$, expansion
\[
 g_Q=\frac1{pz}+\lambda_Q+p\mu_Qz+O(z^2).
\]
The Weierstrass parameter used there agrees with $z$ through degree
three, so these coefficients are unchanged.  Uniqueness gives
$r_a/r_0=pg_{Q_a}$, whence $s_a=p\lambda_{Q_a}$ and
$t_a=p^2\mu_{Q_a}$.  They are precisely the Eisenstein coordinates
of $(C,\omega/p,(Q_a))$.

The fixed Fourier relation selects the required Weil pairing.
Write the pairing of the reconstructed level structure as the $d$th
power of the fixed pairing, for $d\in(\Z/p\Z)^\times$.
Its own Fourier identity and $Fs=-ps$ imply $s_{da}=s_a$.
Equation~\eqref{eq:prime-elimination} then gives $w_{da}=w_a$.
On an elliptic curve equality of Weierstrass coordinates gives
$dQ_a=\pm Q_a$, so $d=\pm1$.  Oddness and $s\ne0$ exclude $d=-1$.
Every noncuspidal complex point of $S$ is consequently on the specified
full-level modular curve.

This proves equality of complex points away from the cusps.  To prove
equality of schemes there, we apply the same construction to first-order
deformations.  Let
$A=\C[\varepsilon]/(\varepsilon^2)$, and take a first-order deformation
of the affine coordinates above a noncuspidal point.
Lemma~\ref{lem:prime-auxiliary} gives a flat auxiliary curve over $A$.
Its special fiber is smooth, so the family is smooth.  The marked
points are disjoint sections.  At $P_\chi$ choose $a\ne0$ and set
$w=r_0/r_a$.  The hyperplane section is reduced on the special fiber,
so $w$ is a local parameter there.  Smoothness identifies the completed
local ring of the family with $A[[w]]$.  The Laurent substitution gives
$w=\chi(a)^{-1}z+O(z^2)$, whose linear coefficient is a unit; hence it
identifies $A[[w]]$ with $A[[z]]$.  In this completion the vector field
is $d/dz$ and is regular.  Its only possible poles were along $r_0=0$.
Faithful flatness of completion therefore implies that it is regular
in the local ring itself.  Its reduction is nonzero everywhere, so it
generates the relative tangent bundle.  The dual differential is a
nowhere vanishing relative differential $\omega_A$.

The character action remains translational.  After composing an
automorphism with the inverse translation of its value at the origin,
it fixes the origin and reduces to the identity.  It is the identity
over $A$, since $H^0(C,T_C(-P_1))=0$.  The marked sections are thus
the full $p$-torsion of the family.  The preceding uniqueness argument
for eigenfunctions works over $A$ as well: their difference is globally
regular, hence in $A$, and a nontrivial character acts on it trivially
only if it is zero.  Integrating $\omega_A$ gives a unique formal
parameter and recovers $s_a,t_a$ from the Laurent coefficients.
Thus the tangent space of the affine cone injects into the deformation
space of an elliptic curve with invariant differential and full level
structure.  This space has dimension two.  Fixing the differential
fixes the scale of a short Weierstrass equation, whose two coefficients
may vary freely, and the level structure lifts uniquely since $p$ is
invertible.  Passing to the projective chart $s_a\ne0$ removes the
scaling direction, so the tangent space of $S$ has dimension at most
one.  The embedded $X(p)$ gives local dimension at least one.
Equality of dimension and embedding dimension proves regularity.

\subsection{The cusps and the saturated ideal}

It remains to check the tangent space at the actual cusps.  Symplectic
changes of indices act transitively on them, so consider infinity.
In the normalization above its coordinates are
\begin{align*}
 s_{(0,j)}=c_j=-\frac{\zeta^j+1}{2(\zeta^j-1)},\quad
 &t_{(0,j)}=\frac1{12},\\
 s_{(r,k)}=B_r=\frac12-\frac rp,\quad
 &
 t_{(r,k)}=\frac{B_r^2-1/12}{2}\quad(r\ne0).
\end{align*}
Write $u=ds$, $v=dt$.  Restriction to the vertical cyclic subgroup
is a tangent vector to the $X_1(p)$ equations.  Their smoothness
\cite[Theorem~4.4]{BGP} and their first Fourier coefficients show that,
after subtracting scaling,
\[
 u_{(0,j)}=A(\zeta^j-\zeta^{-j}),\qquad
 v_{(0,j)}=-A(\zeta^j+\zeta^{-j}).
\]
The mixed cyclic equations linearize to
\begin{equation}\label{eq:prime-cusp-tangent}
 (c_j-B_r)u_{(r,k)}-(c_j+B_r)u_{(r,k+j)}
       +v_{(0,j)}+v_{(r,k)}+v_{(r,k+j)}=0.
\end{equation}
Averaging over $k$ makes $v_{(0,j)}$ independent of $j$, and hence
$A=0$, since $p\ge5$.  For fixed $r\ne0$, expand $u_{(r,k)}$ and
$v_{(r,k)}$ in the characters $\zeta^{hk}$.  Let $U,V$ be their
respective coefficients for a nonzero $h$.  Equation
\eqref{eq:prime-cusp-tangent} gives
\[
 V=\left(B_r-c_j\frac{1-\zeta^{hj}}{1+\zeta^{hj}}\right)U.
\]
If $U\ne0$, the multiplier must be independent of $j$.
This is possible only for $h=\pm1$: after clearing denominators,
compare coefficients in
$(z+1)(1-z^h)=c(z-1)(1+z^h)$ modulo $z^p-1$.
For $h\notin\{0,1,-1\}$ the four exponents $0,1,h,h+1$ are distinct
and give inconsistent coefficients.  The multipliers for $h=1,-1$
are respectively $B_r-1/2$ and $B_r+1/2$.

Write $u_{(r,k)}=\sum_h U_{r,h}\zeta^{hk}$, including the zero vertical
row after subtraction of scaling.  Then
$(Fu)_{(R,T)}=p\sum_r U_{r,-R}\zeta^{-rT}$.
Only rows $R=0,1,-1$ can occur in $Fu$, hence also in $u$.
The equation on $R=0$ gives $U_{r,0}=0$ for every $r$ by Fourier
inversion.  Oddness leaves
\[
 u_{(1,k)}=B\zeta^k+C\zeta^{-k},\qquad
 u_{(-1,k)}=-B\zeta^{-k}-C\zeta^k.
\]
On row $1$ the Fourier equation reads
$p(C\zeta^{-T}-B\zeta^T)=-p(B\zeta^T+C\zeta^{-T})$, so $C=0$.
Thus there is at most one projective tangent direction.  It is present,
being the first coefficient in the actual $q^{1/p}$-expansion.
The local dimension and embedding dimension at the cusp are both one;
its local ring is regular.

We have proved that $S$ has exactly the points of the embedded $X(p)$
and is a regular curve at every one of them.  The closed embedding
$X(p)\hookrightarrow S$ is therefore an equality of schemes.
Let $J$ be the homogeneous ideal in the weight-one coordinate ring
obtained by the elimination \eqref{eq:prime-elimination}, and let
$\mathfrak m$ be the ideal of those coordinates.  Scheme equality
means that every homogeneous polynomial vanishing on the Eisenstein
image belongs to $J:\mathfrak m^\infty$.

Finally take the relation $Q=e$ in the theorem.  Replace its
holomorphic weight-two coordinates by \eqref{eq:prime-elimination}.
The resulting quadratic polynomial vanishes on $X(p)$, and hence is
in $J:\mathfrak m^\infty$.  Restoring the eliminated coordinates
shows that a power of the ideal generated by $x_f$, $f\in V_{p,1}$,
carries $\mathcal L_p(Q)-x_e$ into the weight-two cyclic ideal.
The completed level-one weight-two coordinate is free and does not
occur in this holomorphic relation.  Adjoining the higher-weight
variables does not change the membership just obtained.  This proves
Theorem~\ref{thm:prime-weight-two-saturation}.

\end{document}